\documentclass[12pt]{article}
\usepackage{float}
\usepackage{enumitem}
\usepackage{multirow}
\usepackage{comment}

\usepackage[top=1in, bottom=1in]{geometry} 
\usepackage{graphicx} 
\graphicspath{{./images/}} 
\usepackage{url} 

\usepackage[dvipsnames]{xcolor} 

\definecolor{darkblue}{rgb}{0.0, 0.0, 0.55}
\definecolor{cobalt}{rgb}{0.0, 0.28, 0.67}

\newcommand{\demph}[1]{\emph{{\color{RoyalBlue}#1}}}

\usepackage{tikz} 
\usetikzlibrary{arrows.meta,positioning}

\usepackage{amsmath, amsfonts, amsthm, amssymb} 
\usepackage{mathtools} 

\usepackage{enumitem} 

\makeatletter
\def\saveenum{\xdef\@savedenum{\the\c@enumi\relax}}
\def\resetenum{\global\c@enumi\@savedenum}
\makeatother

\usepackage{soul} 
\usepackage{framed} 
\usepackage{subcaption} 
\usepackage{booktabs} 
\usepackage{arydshln} 

\newcounter{todo}

\makeatletter

\newcommand\listtodoname{List of todos}
\newcommand\listoftodos{%
  \section*{\listtodoname}\@starttoc{tod}}
\makeatother
\newcommand{\jose}[1]{{\color{black}#1}}
\newcommand{\bella}[1]{{\color{black}#1}}

\usepackage[colorlinks=true, allcolors=blue]{hyperref}
 
\usepackage[nameinlink]{cleveref} 

\numberwithin{equation}{section}

\usepackage{aliascnt} 
\theoremstyle{definition}
\newtheorem{theorem}{Theorem}[section]

\newaliascnt{definition}{theorem}
\newtheorem{definition}[definition]{Definition}
\aliascntresetthe{definition}
\crefname{definition}{definition}{definitions}
\Crefname{definition}{Definition}{Definitions}

\newaliascnt{lemma}{theorem}
\newtheorem{lemma}[lemma]{Lemma}
\aliascntresetthe{lemma}
\crefname{lemma}{lemma}{lemmas}
\Crefname{lemma}{Lemma}{Lemmas}

\newaliascnt{notation}{theorem}

\aliascntresetthe{notation}
\crefname{notation}{notation}{notations}
\Crefname{notation}{Notation}{notations}

\newaliascnt{corollary}{theorem}
\newtheorem{corollary}[corollary]{Corollary}
\aliascntresetthe{corollary}
\crefname{corollary}{corollary}{corollaries}
\Crefname{corollary}{Corollary}{Corollaries}

\newaliascnt{proposition}{theorem}
\newtheorem{proposition}[proposition]{Proposition}
\aliascntresetthe{proposition}
\crefname{proposition}{proposition}{propositions}
\Crefname{proposition}{Proposition}{Propositions}

\newaliascnt{problem}{theorem}

\aliascntresetthe{problem}
\crefname{problem}{problem}{problems}
\Crefname{problem}{Problem}{Problems}

\newaliascnt{conjecture}{theorem}

\aliascntresetthe{conjecture}
\crefname{conjecture}{conjecture}{conjectures}
\Crefname{conjecture}{Conjecture}{Conjectures}

\newaliascnt{example}{theorem}
\newtheorem{example}[example]{Example}
\aliascntresetthe{example}
\crefname{example}{example}{examples}
\Crefname{example}{Example}{Examples}

\newaliascnt{remark}{theorem}
\newtheorem{remark}[remark]{Remark}
\aliascntresetthe{remark}
\crefname{remark}{remark}{remarks}
\Crefname{remark}{Remark}{Remarks}

\crefname{mainresult}{main result}{main results}
\Crefname{mainresult}{Main Result}{Main Results}

\usepackage{etoolbox} 

\newcommand{\exampleqedsymbol}{\hfill$\diamond$}      

\AtEndEnvironment{example}{\exampleqedsymbol}

\makeatletter
\let\saved@bibitem\@bibitem
\makeatother

\makeatletter
\let\@bibitem\saved@bibitem
\makeatother

\newcommand{\EDD}{\text{EDDeg}}
\newcommand{\PP}{\mathbb{P}}
\newcommand{\Line}[1]{{\ell^{(#1)}}}
\newcommand{\ellLine}{{\ell}}
\newcommand{\Sing}{\text{Sing}}
\newcommand{\segre}{\sigma}

\newcommand{\Deg}{\operatorname{Deg}}
\DeclareMathOperator{\EDdeg}{EDdeg}
\DeclareMathOperator{\PEDdeg}{EDdeg_{proj}}

\newcommand{\bbeta}{{{\boldsymbol{\beta}}}}
\newcommand{\CC}{{{\mathbb{C}}}}

\newcommand{\Gr}{\operatorname{Gr}}

\newcommand{\LThree}{{\texttt{L}^{\texttt{3}}}}
\newcommand{\HInfinity}{H_\infty}

    \newcommand{\PhiCphi}{\Phi_{\camConf}\circ\phi}

\newcommand{\dimPPImage}{h}
\newcommand{\oneImageAmbient}{\PP^\dimPPImage}
\newcommand{\worldVariety}{Y}
\newcommand{\worldSchubert}{\Lambda}

\newcommand{\camConf}{\mathbf{C}}
\newcommand{\scriptsquare}{%
  \mathbin{\vcenter{\hbox{$\scriptscriptstyle\square$}}}%
}

\newcommand{\snap}[2]{#1 \scriptsquare #2}

\newcommand{\dimPPWorld}{N}
\newcommand{\worldAmbient}{\PP^{\dimPPWorld}}
\newcommand{\imageVariety}{Z}

\newcommand{\plucker}{\iota}
\newcommand{\grCamConf}{\mathbf{D}}

\newcommand{\Aff}{{\operatorname{Aff}}}
\newcommand{\reg}{{\operatorname{reg}}}
\newcommand{\sing}{{\operatorname{sing}}}
\newcommand{\parametricMap}{f}

\newcommand{\field}{\mathbb{K}}
\newcommand{\bfx}{\overline{\mathbf{x}}}
\newcommand{\dimImageVariety}{\dim{\imageVariety}}

\usepackage{tikz-cd}
\usetikzlibrary{arrows.meta}
\newcommand{\hCamConf}{\camConf^{\dimPPImage}}

\newcommand{\hGCamConf}{\grCamConf^\dimPPImage}

\newcommand{\affEDdeg}[1]{\operatorname{affEDdeg}(#1)}

\newcommand{\irreducibleCameraSpace}{\mathcal{A}}

\newcommand{\Xhn}[2]{X_{#1,#2}}

\title{\vspace{-.1in}\jose{The Euclidean distance degree of one-parameter anchored multiview varieties}}

\date{}
\author{Bella Finkel\footnote{This material is based upon work supported by the National Science Foundation Graduate
Research Fellowship Program under Grant No. 2137424, \jose{as well as National Science Foundation Grant No. 2510307}. Any opinions, findings, and conclusions
or recommendations expressed in this material are those of the authors and do not necessarily
reflect the views of the National Science Foundation. Support was also provided by the Graduate
School and the Office of the Vice Chancellor for Research at the University of Wisconsin-Madison
with funding from the Wisconsin Alumni Research Foundation.} 
 \,\,and 
Jose Israel Rodriguez\footnote{This research was partially supported by the Alfred P. Sloan Foundation and by the National Science Foundation Grant No. 2510307.}}
\begin{document}

\maketitle
    \vspace{-24pt}
\begin{abstract}

    Multiview varieties are mathematical models for the set of image feature correspondences that can be produced by a given camera arrangement. They possess an invariant known as their Euclidean distance (ED) degree, which measures the algebraic complexity of determining the 3D features that minimize the reprojection error when computing the scene structure by triangulation.
    In this article, we prove a formula for the ED degree of curves parameterized by rational functions with mild genericity assumptions.  We apply our results to resolve conjectures 
    on one-dimensional line multiview varieties from computer vision proposed by Duff and Rydell. 
\end{abstract}

\section*{Introduction}
One of the most impactful applications of algebraic geometry is in computer vision, leading to the development of the field now known as \emph{algebraic vision}.
A featured problem in this area is to determine the number of critical points \jose{arising in} triangulation and, more generally, \jose{in multiview triangulation} in multiple view geometry (see~\cite{Hartley-Zisserman-book-multiple-view-geometry} for a textbook reference). 
\jose{This amounts to minimizing 
a least squares reprojection error, and is mathematically framed as minimizing the squared Euclidean distance function to an algebraic set known as the (affine) multiview variety.}
Recent work by Duff and Rydell introduced an extensive catalogue of multiview varieties motivated by reconstruction applications 
\cite{Faugeras-version-1,HartleySchaffalitzky,Wolf, LongKanade} and \jose{recent} theoretical advances \cite{breiding2024line,BRST2023-line-view,Ito_Miura_Ueda_2020,rydell2024projectionshigherdimensionalsubspaces,Rydell_2023_ICCV}. In turn, this led to questions and conjectures on the Euclidean distance degrees for these models. 
In this article, we prove two of these conjectures {using a new result on the ED degree of curves in products of projective spaces}.

\begin{figure}[b!]
    \centering
    \includegraphics[width=0.49\linewidth]{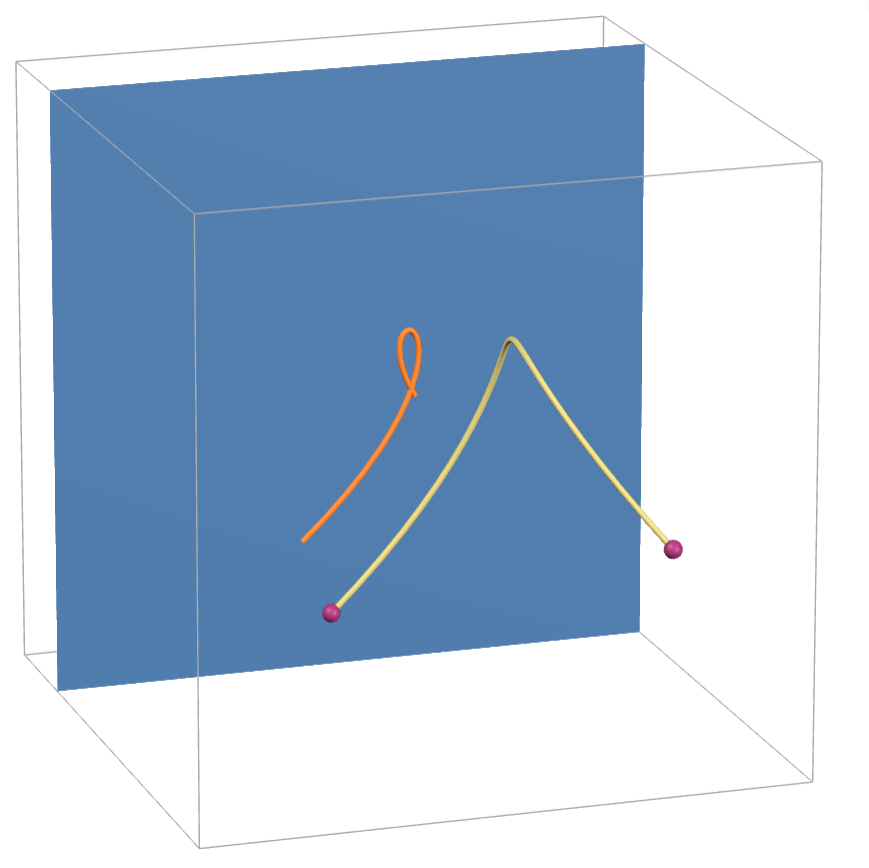}
    \includegraphics[width=0.49\linewidth]{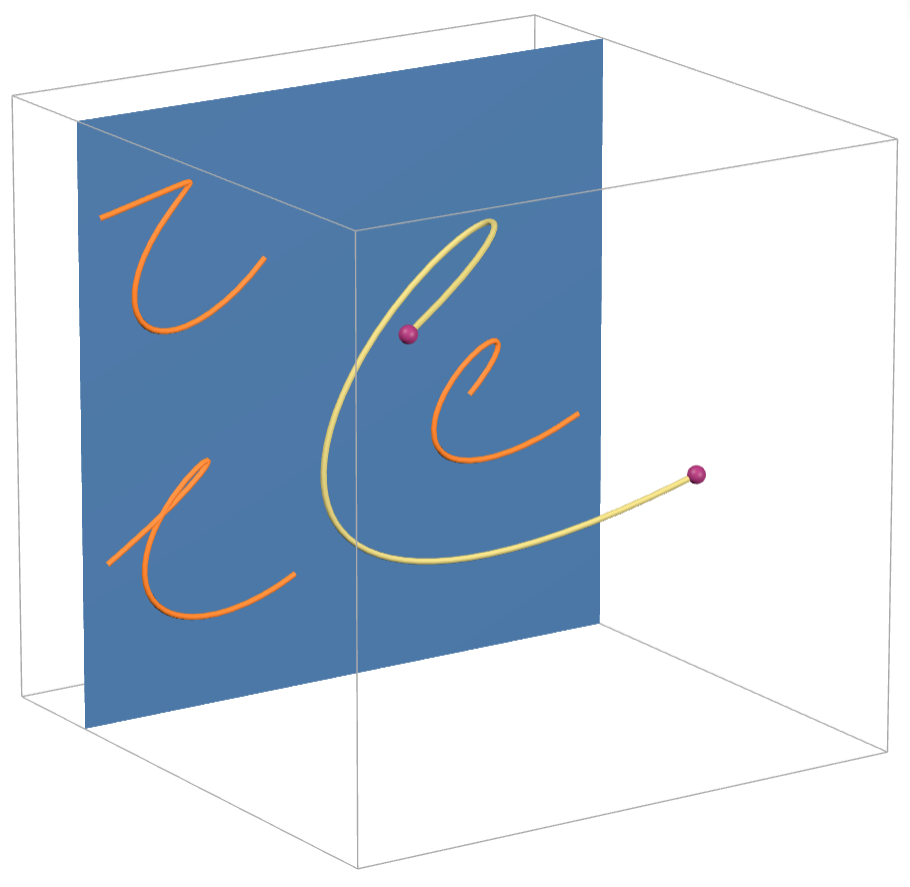}
    \caption{A yellow cubic curve with its orange projection in the blue image~plane (left). 
    A yellow quartic curve with three orange projections in the blue image~plane (right). 
    }
    \label{fig:cubic-quartic-curve}   
\end{figure}

Our paper is structured as follows. In \Cref{s:preliminaries} we provide the necessary preliminaries for the tools we use to prove our main results. We also highlight prior work in the algebraic vision community to put our results into context and address some of the subtleties in the proofs. 
In \Cref{s:multiview-curves} we prove a general statement on the Euclidean distance degree of the {multiview variety of a rational curve of arbitrary degree}.  {In~\Cref{s:ed-degree-conjectures} we have the main application of our result: settling the conjectures
\cite[Conjecture 7.4.5 and Conjecture 7.4.6]{duff-rydel-2024metricmultiviewgeometry} of Duff--Rydell on one-dimensional multiview varieties.}  
\jose{In \Cref{sec:One-parameterFamilies} we apply \Cref{corollary:n-2-is-enough} to determine the ED degree of a class of one-parameter families of lines in the Grassmannian.} 
In \Cref{s:conclusion} we conclude with a brief discussion on future~directions. 

\section{Preliminaries}\label{s:preliminaries}

\subsection{Multiprojective varieties and multidegrees}\label{sss:multiprojective}

\jose{Solution sets of polynomial equations are algebraic varieties. When each polynomial is homogeneous, the solution set is invariant under rescaling, which gives rise to a projective variety. When the defining equations are homogeneous with respect to multiple independent rescalings,
we obtain a subvariety of a product of projective spaces known as a \demph{multiprojective variety}.
In this subsection we review some basics of multiprojective varieties and their multidegrees, and in the next section we see how they arose in computer vision as multiview varieties.
} 

Let $\imageVariety \subset (\mathbb{P}^{\dimPPImage})^n$ be such a variety.
Its defining ideal lies in the multigraded ring
\[
\field[\bfx_1,\dots,\bfx_n]\quad\text{ where }\quad\bfx_i=(x_{i,0},\dots,x_{i,\dimPPImage})\]
and $\field$ is a field.
If $\imageVariety$ is irreducible and of dimension $\dimImageVariety$,
then the multidegree of $\imageVariety$ is the polynomial \[\Deg(\imageVariety):=\sum_{\alpha_1+\cdots+\alpha_n=\dimImageVariety}c_\alpha \jose{{T_1^{\dimPPImage-\alpha_1}\cdots T_n^{\dimPPImage-\alpha_n}}}\]
in the quotient ring  $\mathbb{Z}[T_1,\dots, T_n]/\langle T_1^{\dimPPImage+1},\dots, T_n^{\dimPPImage+1}\rangle$.
\jose{This quotient reflects the geometric fact that one cannot have a nonempty intersection of a factor $\oneImageAmbient$ with more than $\dimPPImage$ general hyperplanes.}
We will use the standard multi-index notation, 
so that if  $\alpha=(\alpha_1,\dots,\alpha_n)$ then $T^\alpha:=T_1^{\alpha_1}\cdots T_n^{\alpha_n}$.
The coefficient $c_\alpha$ is the number of intersection points of $\imageVariety$ with a general linear space of the form 
\[
\mathcal{L}_\alpha:=\cap_{i=1}^n V(\ell_{i,1},\dots, \ell_{i,\alpha_i})
\]
where $\ell_{i,j}$ is a general linear form in the unknowns $\bfx_i$.
One should think of the exponents of the terms of the polynomial as recording how many hyperplanes we impose in each factor, and of the coefficients as counting the size of the corresponding intersection.

If $\imageVariety$ is reducible, then the multidegree of $\imageVariety$ is defined to be the sum of the multidegrees of the irreducible components. 
If $\imageVariety_1$, $\imageVariety_2$ are each equidimensional 
and intersect transversally, 
then 
\[\Deg(\imageVariety_1\cap\imageVariety_2)=\Deg(\imageVariety_1)\cdot \Deg(\imageVariety_2).\]

\begin{example}
{The following shows how multidegrees behave under intersecting bihomogeneous hypersurfaces.}
    Consider the bi-homogeneous polynomials     
\[
    f_d=x_0 y_0^d + x_1 y_1^d+ x_2 y_2^d+x_3 y_3^d,\quad d=0,1,2,3.
\]
Each polynomial $f_d$ defines a hypersurface in $\PP^3\times\PP^3$.
The multidegree of $V(f_d)$ is $T_1+dT_2.$
The multidegree of $V(f_1,f_2,f_3)$ is 
    \[
    (T_1+1T_2)(T_1+2T_2)(T_1+3T_2) =  
        T_1^3+ 6T_1^2 T_2+ 11 T_1T_2^2 + 6T_2^3
    \] 
   and the multidegree of  
    $V(f_1,f_2,f_3)\cap V(f_0)$
    is 
 $6T_1^3 T_2+ 11 T_1^2T_2^2 + 6T_1T_2^3.$
 
 \jose{This example shows how multidegrees behave under intersecting bihomogeneous hypersurfaces.}
\end{example}

\subsection{Multiprojective versus affine multiview varieties}\label{ss:multiprojective}

In multiview geometry, the natural setting is multiprojective space. However, optimization problems such as Euclidean distance minimization are typically formulated in an affine chart. 
\jose{We now describe how multiview varieties appear in both multiprojective and affine settings.
In computer vision, these encode the algebraic constraints satisfied by feature correspondences across multiple images.
}

     An arrangement $\camConf = (C_1,\dots,C_n)$ of full rank $3\times 4$ matrices, called \demph{cameras}, gives a \demph{multiprojective multiview variety of $\PP^3$} 
     as the Zariski closure of the image of the rational map 
\begin{equation}\label{eq:phi-P3-to-P2n}
\begin{aligned}
    \Phi_\camConf:\PP^3 &\dashrightarrow (\PP^2)^n\\
    z &\mapsto (C_1 z, \dots, C_n z).
\end{aligned}
\end{equation}

As in \Cref{sss:multiprojective}, we take the coordinates of $(\PP^2)^n$ to be $\bfx_1,\dots, \bfx_n$.
The multidegree of the image of \eqref{eq:phi-P3-to-P2n} is 
\begin{equation}\label{eq:multidegree-phi-P3-to-P2n}
    \sum_{ \jose{\max\{\alpha_1,\dots,\alpha_n\}\leq 2, }\,\,
    \alpha_1+\cdots+\alpha_n=3
    }
1\cdot {T_1^{2-\alpha_1}\cdot T_2^{2-\alpha_2}\cdot\cdots\cdot T_n^{2-\alpha_n}},\quad n>1.
\end{equation}
The exponent $2-\alpha_i$ indicates that $\alpha_i$ hyperplanes were used in the $i$th factor.
Ideals and multidegrees of multiview varieties are also well studied in algebraic vision, and for more details we highlight~\cite{AST-multidegree-hilbert-scheme,BRST2023-line-view, EscobarKnutson2017-multidegree}.
For a textbook introduction to multidegrees, see~\cite[Chapter 8]{Miller-Sturmfels-combinatorial}.  

\newcommand{\dotPoduct}[2]{#1\cdot#2}
The \demph{affine multiview variety of $\PP^3$} is obtained by restricting each factor of $(\PP^2)^n$ to an affine chart. 
Our convention is to take the charts given by $x_{i,0}=1$ for $i=1,\dots,n$. Under this choice, the affine multiview variety has the following parameterization:
\[
z\mapsto 
        \begin{bmatrix}
            \displaystyle{\frac{
                \dotPoduct{ C_1^{(1)} }{z}
                }{\dotPoduct{ C_1^{(0)}}{z}}}
                &\dots &
                            \displaystyle{
                            \frac{
                \dotPoduct{ C_n^{(1)} }{z}
                }{\dotPoduct{ C_n^{(0)}}{z}}}\\ \\
            \displaystyle{\frac{
                \dotPoduct{ C_1^{(2)} }{z}
                }{\dotPoduct{ C_1^{(0)}}{z}}
                }
            &\dots &
            \displaystyle{\frac{
                \dotPoduct{ C_n^{(2)} }{z}
                }{\dotPoduct{ C_n^{(0)}}{z}}}
        \end{bmatrix}
        \in \field^{2\times n}.
\]
In general, if $\imageVariety$ is a multiprojective subvariety
in $(\oneImageAmbient)^n$ 
then
we denote its restriction to the conventional affine chart 
by $\imageVariety_{\Aff}\subset (\field^{\dimPPImage})^n$. 
\jose{This corresponds to working with inhomogeneous image coordinates, where reprojection error is naturally defined.}

\subsection{Euclidean distance degree}\

\newcommand{\XR}{X_\mathbb{R}}

The \demph{Euclidean distance degree} (ED degree) of an algebraic variety measures the algebraic complexity of the problem of minimizing Euclidean distance to that variety. 
We recall its definition in this subsection.

If $X$ is an affine variety in $\CC^\dimPPImage$, then we define 
$\XR:=X\cap \mathbb{R}^\dimPPImage$
to be the real points of $X$.
We call 
$u\in\mathbb{R}^\dimPPImage$  a data point. 
The optimization problem 
\[
\min_{x\in \XR} \Vert x-u\Vert^2
\]
has first-order optimality conditions. 
Namely, for a smooth point $x\in \XR$, we require
$(x-u)^Tv= 0$ for all $v$ in the tangent space $T_x\XR$ of $\XR$ at $x$. 
This condition gives rise to a polynomial system \cite[Equation 2.1]{DHOST2016-ed-degree} that has $x$ as the unknown and $u$ as the given generic data. 
The number of complex (real or nonreal) solutions to this polynomial system is called the Euclidean distance degree of $\XR$~\cite{DHOST2016-ed-degree}. 
The ED degree is a measure of complexity for minimizing the squared Euclidean distance from $\XR$ to $u$. 

\begin{figure}[hbt!]
    \centering
    \includegraphics[width=0.5\linewidth]{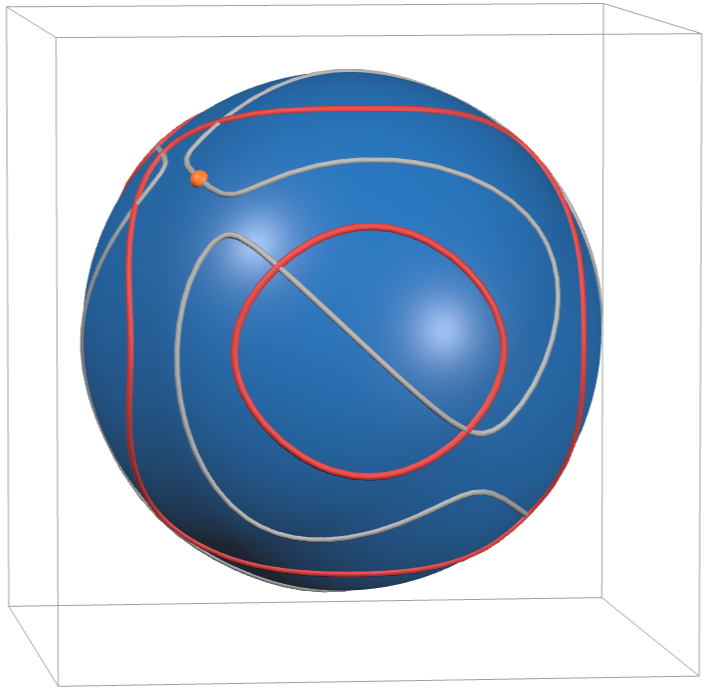}
    \caption{ 
    The red curve is $\XR$ and the real ED critical points for the distance function to the orange point is the set of intersection points with the silver curve. 
    } 
    \label{fig:curve-on-sphere}
\end{figure}

\begin{example}
     Let $X$ denote the variety of the ideal generated by
     \[
     g_1=x^2+y^2+z^2-1,\quad g_2=x^{4}+y^{4}+z^{4}+ax^{2}y^{2}+by^{2}z^{2}+cz^{2}x^{2}-s\left(x^{2}+y^{2}+z^{2}\right)^{2}\]
     where $(a,b,c,s)=(-0.8, -0.74, -0.5,0.32)$. 
     Then, for a generic data point $u$, the ideal of the critical points is generated by
     $g_1,g_2$ and the determinant \jose{of the $3\times3$ matrix}
     \[\begin{bmatrix}
      \nabla d_u &\nabla g_1 & \nabla g_2   
     \end{bmatrix}\]
     where $d_u$ is the squared Euclidean distance function from the data $u$ and $\nabla g$ denotes the $3\times 1$ matrix of partial derivatives of $g$.
        In \Cref{fig:curve-on-sphere}, 
        the silver curve is the restriction of the variety of the determinant to the sphere and the real ED critical points are the intersection points of the silver and red curves.

     Note, when $X$ is singular an additional step of saturation needs to be performed to obtain the ideal of critical points. 
\end{example}

\medskip
Beyond its algebraic definition, the Euclidean distance degree also has a topological interpretation. For more details, we refer to the survey article~\cite{MRW-survey}. 

\newcommand{\Eu}{\operatorname{Eu}}
\begin{theorem}[\cite{MRW-multiview}]\label{theorem:ed-affine}
    Suppose $X$ is an irreducible  closed subvariety of $\mathbb{C}^{\dimPPImage}$. Then for general $\boldsymbol{\beta}:=(\beta_0,\dots,\beta_\dimPPImage)\in \mathbb{C}^{\dimPPImage+1}$ we have
    \[\EDD(X) = 
        (-1)^{\dim X}
        \chi(\Eu_X|_{X\cap U_{\boldsymbol{\beta}}})
        \]
    where
    $\Eu_X$ is the local Euler obstruction function on $X$ and
    \[ U_{\boldsymbol{\beta}}
        :=\mathbb{C}^\dimPPImage
            \setminus\{
    (z_1,\dots,z_\dimPPImage)\in
        \mathbb{C}^\dimPPImage: 
            \sum_{i=1}^\dimPPImage(z_i-\beta_i)^2+\beta_0=0
    \}.\]
    
\end{theorem}

For many of our purposes, $X$ is  smooth and it suffices to consider the
 Euler--Poincar\'e characteristic $\chi$.

\begin{theorem}[\cite{MRW-multiview}]\label{thm:MRW-smooth}
    Suppose $X$ is a smooth closed subvariety of $\mathbb{C}^\dimPPImage$. 
    Then for general $\boldsymbol{\beta}:=(\beta_0,\dots,\beta_\dimPPImage)\in \mathbb{C}^{\dimPPImage+1}$ we have
    \[\EDD(X) = 
        (-1)^{\dim X}
        \chi(X\cap U_{\boldsymbol{\beta}})
        \]
    where $U_{\boldsymbol{\beta}}:=\mathbb{C}^\dimPPImage\setminus\{
    (z_1,\dots,z_\dimPPImage)\in \mathbb{C}^\dimPPImage: \sum_{i=1}^\dimPPImage(z_i-\beta_i)^2+\beta_0=0
    \}$.
\end{theorem}

{We apply these results to compute the ED degrees of one-dimensional multiview varieties in} \Cref{s:multiview-curves} and \jose{resolve computer vision conjectures in \Cref{s:ed-degree-conjectures}.} 

\subsection{The projective ED degree}\label{sss:projective-ed}

The \demph{projective ED degree} provides a natural way to study ED degrees of affine cones by avoiding the singularity at the cone point and enabling the use of intersection-theoretic tools.
In this subsection, we recall this invariant and prior work for curves. 
\jose{Although our main results concern affine ED degrees, the projective ED degree also plays an important role because many multiview varieties arise as affine cones or have natural projective closures.}

\smallskip
Let $\imageVariety\subset \oneImageAmbient$ be an irreducible closed subvariety of codimension $c$ and let $I_\imageVariety$ be its defining ideal.
If the homogeneous polynomials $g_1,\dots,g_s\in \mathbb{C}[z_0,z_1,\dots,z_\dimPPImage]$ generate $I_\imageVariety$, then the \demph{projective Euclidean distance correspondence} is a biprojective variety in $\PP^{\dimPPImage}\times\PP^{\dimPPImage}$ defined by the saturated ideal
\begin{equation}
\label{eq:critideal2}
\biggl( I_\imageVariety + \biggl\langle
\hbox{$(c+2) \times (c+2) $-minors of $ \begin{pmatrix} u \\ z \\ J(g) \end{pmatrix} $} \biggr\rangle
\biggr) : \bigl(I_{\imageVariety_{\rm sing}} \cdot  \langle z_0^2 + \cdots + z_\dimPPImage^2 \rangle \,\bigr)^\infty 
\end{equation}
 where 
 $J(g)$ is the $s\times (\dimPPImage+1)$ matrix of partial derivatives of $g$ and
 the singular
locus $\imageVariety_\sing$ of $\imageVariety$ is given by the ideal
$$ I_{\imageVariety_\sing} \,\,= \,\,
I_\imageVariety + \bigl\langle \hbox{$(c \times c)$-minors of  $J(g)$} \bigr\rangle.$$

The \demph{ED bidegree} is defined to be multidegree of the ED correspondence. 
The coefficient of  $T_2^{\dimPPImage}$ is the
(projective) ED degree of $X$, and is denoted by $\PEDdeg(X)$.

The affine cone over $\imageVariety$ is
\[
\{(z_0,\dots,z_\dimPPImage)\in\mathbb{C}^{\dimPPImage+1}: [z_0:\dots:z_\dimPPImage]\in \imageVariety\}.
\]
The ED degree of this affine cone agrees with $\PEDdeg(\imageVariety)$ \cite{DHOST2016-ed-degree}. The advantage of working projectively is to avoid the singularity at the cone point. It also allows the use of several topological formulas. 

\begin{theorem}\cite[Theorem 1.3]{MRW2021-edprojective}
    \label{MRW:ed-projective}
Let $\imageVariety$ be any irreducible closed subvariety of $\oneImageAmbient$. Then for a general $\bbeta=(\beta_0,\ldots,\beta_{\dimPPImage})\in \CC^{{\dimPPImage+1}}\setminus\{0\}$ we have
\[\label{e4i} 
\PEDdeg(\imageVariety)=(-1)^{\dim X}\chi\big({\rm Eu}_{\imageVariety}|_{\mathcal{U}_\bbeta}\big),
\]
where 
    \[ \mathcal{U}_\bbeta:=\PP^{{\dimPPImage}} \setminus (Q \cup H_{\bbeta})
    \]
with $Q$ denoting the 
quadric $\{[z_0:\dots:z_{\dimPPImage}]\in \PP^{{\dimPPImage}}: \sum_{i=0}^{\dimPPImage} z_i^2=0\}$ and 
\[H_{\bbeta}:= \{[z_0:\dots:z_\dimPPImage]\in \PP^{{\dimPPImage}}:\beta_0 x_0+\ldots + \beta_{{\dimPPImage}}x_{\dimPPImage}=0 \}.\]
\end{theorem}

 When $\imageVariety\subset\oneImageAmbient$ is smooth, the formula simplifies.

\begin{theorem}\cite[Theorem 8.1]{aluffi-harris-projective-ed} \label{theorem:smooth-projective}
Let $\imageVariety$ be a smooth subvariety of $\oneImageAmbient$, and assume that $\imageVariety \nsubseteq Q$, where 
$Q=\{[z_0:\ldots:z_\dimPPImage]\in \oneImageAmbient \mid z_0^2+\ldots + z_\dimPPImage^2=0\}$ is the 
quadric in $\oneImageAmbient$. 
Then
\[\label{e3i}
\PEDdeg(\imageVariety)=(-1)^{\dim \imageVariety}\chi\big(\imageVariety \setminus (Q \cup H)\big),
\]
where $H$ is a general hyperplane. 
\end{theorem}

As noted in \cite[Section 9]{aluffi-harris-projective-ed}, an immediate application of the previous theorems gives a formula for the ED degree of a smooth curve $\imageVariety$ in $\oneImageAmbient$ of degree $d$:
\begin{equation}\label{eq:ed-smooth-curve}
\PEDdeg(\imageVariety)=d+\#(Q\cap \imageVariety)-\chi(\imageVariety).
\end{equation}

The next examples illustrate some of the discrepancies that occur between an affine variety's ED degree (which is our setting) and the ED degree of related projective varieties.
\begin{example}[Standard grading giving discrepancy]\label{ex:cubic-10-14}
    Consider the curve $Z$ in $\PP^2\times \PP^2$ parameterized by the map
    \[
        \PP^1\to Z\subset \PP^2\times\PP^2,
        \]
    \[
    [s:t]\mapsto \left(\begin{bmatrix}
        s^3+2t^3\\
        s^2t+3t^3\\
        st^2+5t^3
    \end{bmatrix},
    \begin{bmatrix}
        s^3+7t^3\\
        s^2t+11t^3\\
        st^2+13t^3
    \end{bmatrix}
    \right).
    \]
    The affine variety $Z_\Aff\subset\field^4$ has ED degree  10. On the other hand, we can regard the ideal of  $Z$ as a homogeneous ideal in the usual $\mathbb{Z}$-grading. This defines a projective surface in $\PP^5$ that has ED degree 14. 
\end{example}

\jose{The ED degree is not invariant under projective closure. Even for classical multiview varieties, the affine and projective versions may have drastically different ED degrees.}

\begin{example}[Projective closure discrepancy]\label{ex:projective-closure-discrepancy}
    The affine $n$-view variety discussed in \Cref{ss:multiprojective}
    is known
        \cite{MRW-multiview} 
    to have ED degree 
     \begin{equation}\label{eq:big-conjecture}
     \frac{9}{2}n^3-\frac{21}{2}n^2+8n-4.
     \end{equation}
     \jose{This cubic polynomial in $n$ counts the number of (real or nonreal) critical points of the distance function from a general data point to an $n$-view variety when the cameras are in generic position. The values of this polynomial for small $n$ was determined by \cite{SSN2005-three-view}. In particular, for $n=3$, there are forty seven critical points, and \cite{SSN2005-three-view} showed how to minimize the distance function to a three-view variety by solving a $47\times47$ eigenvalue problem.} 
     
     \jose{In \cite{DHOST2016-ed-degree}, it was noticed that the projective closure in $\PP^{5}$ of a three-view variety has ED degree 56. 
     In \cite[Theorem~1]{Harris03062018}, they found the projective closure of the $n$-view variety in $\mathbb{P}^{2n}$ has  ED degree}  
     \[
     6n^3-15n^2+11n-4.\vspace{-24pt}
     \]
     
\end{example}

\subsection{More multiview varieties and triangulation}

\jose{Anchored multiview varieties arise naturally when the world contains known structures, such as a curve, line, or surface. These models describe all possible image observations of points constrained to lie on $\worldVariety$.}
In this subsection, we recall the definition of \demph{anchored point multiview varieties}. 
To stay organized and minimize notation, we focus on the camera arrangement and the image it induces.

\begin{definition}\label{def:MultiviewVariety-V2}
    Let $\camConf = (C_1,\dots,C_n)$ be a \demph{camera arrangement}, where each $C_i$ is a full-rank $(\dimPPImage+1)\times(\dimPPWorld+1)$ matrix.
    Define the rational map 
    \begin{equation}\label{eq:PhiCam-V2}
        \begin{aligned}
        \Phi_\camConf : 
            \worldAmbient \dashrightarrow (\oneImageAmbient)^n, \quad X \mapsto (C_1 X,\dots,C_n X).
        \end{aligned}
    \end{equation}
     The \demph{(point) multiview variety of $\worldAmbient$} with respect to $\camConf$ is 
        \[ 
        \snap{\camConf}{\worldAmbient} 
            := 
        \overline{\operatorname{im}(\Phi_\camConf)} \subset (\oneImageAmbient)^n,
        \] 
        where the overline denotes the Zariski closure. 
    If $\worldVariety \subseteq \worldAmbient$ is a subvariety, 
         define 
    \[ 
    \snap{\camConf}{\worldVariety} 
        := 
        \overline{\Phi_\camConf(\worldVariety)}
    \]        
    and call this the \demph{(point) multiview variety anchored at $\worldVariety$} with respect to $\camConf$.
     We denote its restriction to the standard affine chart by $(\snap{\camConf}{\worldVariety})_{\Aff}$ and refer to it as the \demph{affine multiview variety of $\snap{\camConf}{\worldVariety}$}. 
\end{definition}

When $\dimPPWorld=3$ and $\dimPPImage=2$, the map $\Phi_\camConf$ models the procedure of taking $n$ pictures with $n$ possibly distinct cameras. A world point $P\in\PP^3$ is mapped to the tuple of images $(C_1P,\dots,C_nP)$ provided that $P$ is not in the kernel of any $C_i$. Points $P\in\PP^\dimPPWorld$ that lie in the kernel of some camera are referred to as being at a \demph{camera center}, where $\Phi_\camConf$ is undefined. {That is, the camera centers constitute the base locus of $\Phi_\camConf$.}

\begin{example}[A one-dimensional multiview variety]\label{example:three-n-minus-two}
    If $\worldVariety$ is a line in $\PP^3$, consider the anchored multiview variety
    $\snap{(C_1,\dots,C_n)}{Y}$.
    The affine multiview variety of $\snap{(C_1,\dots,C_n)}{Y}$ has ED degree $3n-2$. 
    This result appears in \cite[Theorem 1.7]{Rydell_2023_ICCV} and 
    serves as a baseline case for understanding the complexity of triangulating points on curves.
    This is the degree one case of \Cref{thm:ed-curve}.
\end{example}

In~\cite{duff-rydel-2024metricmultiviewgeometry}, they avoid the subscript $\operatorname{Aff}$ by \emph{defining} the ED degree of 
$\snap{\camConf}{\worldAmbient}$
to be
the ED degree of 
$(\snap{\camConf}{\worldAmbient})_{\operatorname{Aff}}$.
They mention how this choice is inconsistent with \cite{DHOST2016-ed-degree} in the one view setting. 
To avoid confusion,
we opt to use the following notation. If $X$ is a multiprojective variety, 
then 
    \begin{equation}
        \affEDdeg{X} :=  \EDdeg(X_\Aff)
    \end{equation}
where $X_\Aff$ is the restriction of $X$ to the  affine chart described in \Cref{ss:multiprojective}.

\subsection{A Schubert variety in $\Gr(1,\PP^3)$}\label{ss:schubert}

In this subsection, we \jose{consider} a subvariety of the Grassmannian that will induce the multiprojective varieties appearing in our results in \Cref{s:ed-degree-conjectures}.
Let $\Gr(1,\PP^3)$ denote the Grassmannian of lines in $\PP^3$. Given three skew lines $\Line{1},\Line{2},\Line{3}$ we \jose{have}

\begin{equation}\label{eq:LThree-defined}
    \LThree(\Line{1},\Line{2},\Line{3}):=\{
        \ellLine \in\Gr(1,\PP^3): \ellLine\cap \Line{i}\neq \emptyset,\, i=1,2,3 
        \}.
\end{equation}
This is an example of a Schubert variety~\cite{KL1972-Schubert}. We often abbreviate $\LThree(\Line{1},\Line{2},\Line{3})$
by $\LThree$. 
To understand the geometry of $\LThree$, we recall two classical results about lines and quadrics in $\PP^3.$

\newcommand{\LQuadric}{S}

 \begin{theorem}\label{theorem:containing-LThree}
    There is a unique smooth quadric $\LQuadric$ that contains any three given pairwise disjoint lines in $\mathbb{P}^3$. 
\end{theorem}
A consequence of \Cref{theorem:containing-LThree} is that any line $\ellLine\in\Gr(1,\PP^3)$ satisfying the incidence condition~\eqref{eq:LThree-defined} must intersect the quadric at three points and is therefore contained in it. 
Thus, 
$\LThree$
 lies within the family of lines contained in the quadric.
The next theorem characterizes the lines contained in a quadric. 

\begin{theorem}\label{theorem:LThree-contains}
    A smooth quadric $\LQuadric$ in $\PP^3$ contains two {rulings}.
    Moreover, 
    any two lines from different 
    {rulings}
    intersect in exactly one point, and 
    any two lines from the same 
    {ruling}
    are disjoint.
\end{theorem}
\Cref{theorem:LThree-contains}
implies that the three skew lines in our incidence condition are in the same ruling.  All together, these theorems imply that $\LThree$ is the other ruling.  

The Pl\"ucker embedding 
\[\plucker:\Gr(1,\PP^3)\to \PP^5\] provides us a way to embed $\LThree\subset\Gr(1,\PP^3)$ into $\PP^5$. 
It turns out $\plucker(\LThree)$ is a conic curve in $\PP^5$, and the next example illustrates how to parameterize it.

\begin{example}[Parameterizing $\LThree$]\label{ex:embedding-L3}
In this example, we describe the parametrization of $\plucker(\LThree)$
where $\Line{i}$ is the line 
\begin{equation}\label{eq:chose-three-skew-lines}
\PP^1\to \PP^3,\,\quad
[s:t]\mapsto [su_i:sv_i:tu_i:tv_i]\quad\text{ with fixed }\quad  [u_i:v_i]\in \PP^1.
\end{equation}
The unique smooth quadric $S$ containing these lines is 
the image of the Segre embedding
\begin{align*}
    \segre:\PP^1\times\PP^1&\to\PP^3\\
    ([s:t],[u:v])&\mapsto[su:sv:tu:tv]
\end{align*}
so that $S$ is defined by the equation $x_0x_3-x_1x_2=0$.
The line $\Line{i}$ is spanned by the points $[u_i:v_i:0:0]$ and $[0:0:u_i:v_i]$, and
we get a ruling in $S$ containing the three skew lines
by varying the point $[u:v]\in \PP^1$. 
We get the second ruling, $\LThree$, by varying $[s:t]\in \PP^1$ to get lines spanned by 
the points $[s:0:t:0]$ and $[0:s:0:t]$.

In terms of Pl\"ucker coordinates on $\PP^5$, lines in $\iota(\LThree)$ are given by
\begin{equation}\label{eq:L3Coordinates-V2}
    \begin{bmatrix}
    p_{12}\\
    p_{13}\\
    p_{23}\\
    p_{14}\\
    p_{24}\\
    p_{34}
\end{bmatrix}=\begin{bmatrix}
    s^2\\ 0\\ -st\\ st\\  0\\t^2
\end{bmatrix}
\end{equation}

\smallskip
To parameterize the lines in $\plucker(\LThree)$ by $\PP^1$, we use a composition of the Veronese embedding 
\[\nu([s:t])=\begin{bmatrix}
    s^2\\ st \\ t^2
\end{bmatrix}\]
with a linear map, \jose{denoted by $B_\LThree$},
so that the parameterization map $\parametricMap:\PP^1\to \plucker(\LThree)$ takes $[s:t]$ to 
    $B_{\LThree}\cdot \nu([s:t])$.

Note, our choice for the three skew lines $\Line{1},\Line{2},\Line{3}$ in \eqref{eq:chose-three-skew-lines} was made without loss of generality as the unique quadric containing any other choice of three skew lines is projectively equivalent to $V(x_0x_3-x_1x_2)$. 
\end{example}

\begin{remark}
If $\grCamConf$ is an arrangement of cameras of size $(H+1)\times 6$, 
then $\snap{\grCamConf}{\plucker(\LThree)}\subset (\PP^H)^n$ is a point multiview variety anchored at $\plucker(\LThree)$. 
This variety was previously studied in \cite{duff-rydel-2024metricmultiviewgeometry} when
$\grCamConf$ is an arrangement of so-called \demph{wedge cameras} (\Cref{definition:wedge-cameras}).
We revisit this construction and show how it relates to \emph{line multiview varieties} 
in \Cref{s:ed-degree-conjectures} for $H=\binom{2+1}{2}-1$ and $H=\binom{3+1}{2}-1$.
\end{remark}

\newcommand{\phiWorldVariety}{\parametricMap(\PP^1)}

\section{ED degrees of multiview varieties anchored at rational curves}\label{s:multiview-curves}

\newcommand{\ee}{E}
In this section, 
we let $\worldVariety$ be a degree $\ee$ rational curve in $\worldAmbient$ with $\dimPPWorld\geq 3$ parametrized~by 
\begin{align}\label{eq:phi-P1-to-Pell}
    \parametricMap:\PP^1&\to\worldAmbient\\
    [s:t]&\mapsto \parametricMap(s,t)=[\parametricMap_0(s,t):\parametricMap_1(s,t):\dots:\parametricMap_{\dimPPWorld}(s,t)]\notag.
\end{align}
{We say $\worldVariety$ is the \demph{world variety}.}
We always assume the $\parametricMap_i$ are homogeneous degree $\ee$ polynomials such that $\parametricMap$ has no base locus and
is generically one-to-one.
Unless otherwise stated, we also assume $\worldVariety=\phiWorldVariety$ is smooth or only has nodes as singularities.

As a special case of \Cref{def:MultiviewVariety-V2}, for $\worldVariety$ as in \eqref{eq:phi-P1-to-Pell} 
we denote its multiview variety with respect to the camera arrangement $\camConf$
as $\snap{\camConf}{\phiWorldVariety}$. 
If $\worldVariety$ does not contain any points in a camera center of the arrangement, then the image $\Phi_\camConf(\phiWorldVariety)$ is closed.

\begin{lemma}\label{lemma:smooth-or-nodal}
    {Assume we are given a generic camera arrangement $\camConf=(C_1,\dots,C_n)$ of full rank $(\dimPPImage+1)\times(\dimPPWorld+1)$ matrices $C_i$ with 
    $\dimPPWorld\geq3$ and
    $\dimPPImage\geq2$.}
{Let $\parametricMap:\PP^1\to\worldAmbient$ be as in \eqref{eq:phi-P1-to-Pell} so that $\worldVariety=\phiWorldVariety$ is degree $\ee$ and $\worldVariety$ is either smooth or has only nodal singularities.}
Then  $\snap{\camConf}{\phiWorldVariety}$ is either a smooth curve or only has nodes as singularities. 
\end{lemma}
\begin{proof}
For $n=1$, this is a standard algebraic geometry fact. For instance, 
see \cite[Exercise~3.34]{eisenbud20163264}.
The case when $n>1$ immediately follows.
\end{proof}

\begin{example}[Single view varieties of cubic curves]\label{ex:cubic-curves}
     Let $\worldVariety=\phiWorldVariety$ be the rational cubic curve in $\PP^{3}$ of degree $3$
    given by 
    \[
        [s:t]\overset{\parametricMap}{\mapsto} [t^3:st^{2}:s^2t:s^3].
        \]
    Let $C$ be a general $(\dimPPImage+1)\times {4}$ camera matrix and take the camera configuration to be this single camera ($n=1$).

    There are three different cases: 
    \begin{enumerate}[nosep]
    \item     If $\dimPPImage=1$, then $\snap{\camConf}{\phiWorldVariety}$ is $\PP^1$. The ED degree of $(\snap{\camConf}{\phiWorldVariety})_{\operatorname{Aff}}$ is one. 

    \item If $\dimPPImage=2$, then $\snap{\camConf}{\phiWorldVariety}$ is a  curve in $\PP^2$ with a nodal singularity.  

    \item If $\dimPPImage\geq3$, then $\snap{\camConf}{\phiWorldVariety}$ is a smooth curve in $\oneImageAmbient$. 
    \end{enumerate}
    If $\dimPPImage\geq 2$ then the 
    ED degree of $(\snap{\camConf}{\phiWorldVariety})_\Aff$ is seven. 
    {In \Cref{fig:cubic-quartic-curve} (left) we plot $\snap{\camConf}{\phiWorldVariety}$ in orange 
    for a camera of size $3\times4$ and the curve  $\worldVariety$ in yellow.
    } 
\end{example}

\begin{theorem}\label{thm:ed-curve}
{Assume we are given a generic camera arrangement $\camConf=(C_1,\dots,C_n)$ of full rank $(\dimPPImage+1)\times(\dimPPWorld+1)$ matrices $C_i$ with 
$\dimPPWorld\geq3$ and
$\dimPPImage\geq2$.}
{Let $\parametricMap:\PP^1\to\worldAmbient$ be as in \eqref{eq:phi-P1-to-Pell} so that $\worldVariety=\phiWorldVariety$ is degree $\ee$ and $\worldVariety$ is smooth or has only nodal singularities.}
    Then the ED degree of an affine patch $(\snap{\camConf}{\phiWorldVariety})_\Aff$ of $\snap{\camConf}{\phiWorldVariety}$ is
     \begin{equation}
        \affEDdeg{\snap{\camConf}{\phiWorldVariety}} = 3\ee n-2.
    \end{equation}
\end{theorem}
\begin{proof}

Recall from \Cref{def:MultiviewVariety-V2}
   that a camera $C$ is a full rank $(\dimPPImage+1)\times(\dimPPWorld+1)$ matrix. We 
    write $C^{(i)}$ for the $i$th row of $C$. 
    The subvariety    
    $X:=\snap{\camConf}{\parametricMap(\PP^1)}$
    has the parameterization
    \begin{equation}\label{eq:CamImageProj}
    \Phi_{\mathbf{C}}\circ \parametricMap([s:t])=
        \begin{bmatrix}
            \begin{pmatrix}
                \dotPoduct{C_1^{(0)}}{f}\\
            \dotPoduct{C_1^{(1)}}{f}\\
            \vdots\\
            \dotPoduct{ C_1^{(\dimPPImage)}}{f}
            \end{pmatrix}:\dots:\begin{pmatrix}
                \dotPoduct{C_n^{(0)}}{f}\\
            \dotPoduct{C_n^{(1)}}{f}\\
            \vdots\\
            \dotPoduct{C_n^{(\dimPPImage)}}{f}
            \end{pmatrix}
        \end{bmatrix}\in (\oneImageAmbient)^n.
    \end{equation}
     Let $X_\Aff:=(\snap{\camConf}{\worldVariety})_\Aff$ denote the restriction of $\operatorname{im}\Phi_\camConf$ to an affine patch  with coordinates
    \begin{equation}\label{eq:X-aff-param}
        \begin{pmatrix}
            \begin{pmatrix}
                \frac{\dotPoduct{C_1^{(1)}}{f}}{\dotPoduct{C_1^{(0)}}{f}}\\
                \vdots\\
                \frac{\dotPoduct{C_1^{(\dimPPImage)}}{f}}{\dotPoduct{C_1^{(0)}}{f}}
            \end{pmatrix}, \dots, \begin{pmatrix}
                \frac{\dotPoduct{C_n^{(1)}}{f}}{\dotPoduct{C_n^{(0)}}{f}}\\
                \vdots\\
                \frac{\dotPoduct{C_n^{(\dimPPImage)}}{f}}{\dotPoduct{C_n^{(0)}}{f}}
            \end{pmatrix}
        \end{pmatrix}\in (\field^{\dimPPImage})^n.
    \end{equation}

    In $(\oneImageAmbient)^n$, each factor of $\oneImageAmbient$ gives a hyperplane at infinity. We denote the union of these hyperplanes as $\HInfinity$. If the $i$th factor in $\oneImageAmbient$ has coordinates $\bfx_i=(x_{i,0},x_{i,1},\dots,x_{i,\dimPPImage})$ then
    \[
        \HInfinity=\cup_{i=1}^n V(x_{i,0}).
    \]
    Using the parameterization of $X$ in \eqref{eq:CamImageProj}, 
    we have 
    \[
    X\cap\HInfinity=
    \PhiCphi(S_\infty)
    \text{\quad where\quad}
    S_\infty:= \bigcup_{i=1}^nV( \dotPoduct{C_i^{(0)}}{f} )\subset \PP^1.
    \] 
    If each $\dotPoduct{C_i^{(0)}}{f}$ has distinct roots and 
    \begin{equation}\label{eq:generic-condition-infinity}
        V(\dotPoduct{C_i^{(0)}}{f},\dotPoduct{C_j^{(0)}}{f})=\emptyset
    \end{equation}
    for all $i,j\in[n],\ i\neq j$, 
    then $S_\infty$ has $\ee n$ points.  
    {This holds by genericity of the camera arrangement: the degree $\ee$
    polynomial $\dotPoduct{C_i^{(0)}}{f}$ has a multiple root if and only if its $X$-discriminant vanishes, and $\dotPoduct{C_i^{(0)}}{f},\dotPoduct{C_j^{(0)}}{f}$ share a root if and only if their resultant is nonvanishing \cite[Section 1.1]{GKZ}. 
    }    
    Therefore, if $\parametricMap$ is a generically one-to-one map then 
    $\HInfinity\cap X$ has precisely $\ee n$ points. 
    \medskip

    For 
    $\beta=(\beta_0,\beta_{1,1},\dots,
    \beta_{n,\dimPPImage})\in\mathbb{R}^{n\dimPPImage+1}$, define 
    the hypersurface
   \begin{equation}\label{eq:Q-beta-U-chart}
        Q_{\beta}:=
\left\{
    (z_{1,1},\dots,
    z_{n,\dimPPImage})\in \mathbb{C}^{n\dimPPImage}:
    \sum_{i=1}^n \sum_{j=1}^\dimPPImage(z_{i,j}-\beta_{i,j})^2+\beta_0=0
    \right\}
    \end{equation}
    and the Zariski open set
   \begin{equation}\label{eq:U-beta-U-chart}
        U_{\beta}:=
    \mathbb{C}^{\dimPPImage n}\setminus
    Q_\beta.
    \end{equation}
   Let $S_Q$ be the  set of points in $\PP^1\setminus S_\infty$ 
   that are the zeros of 
    \begin{align*}
    d_\beta = 
    \sum_{i=1}^n \sum_{j=1}^\dimPPImage
        \left(
            \frac{\dotPoduct{C_i^{(j)}}{f}}{\dotPoduct{C_i^{(0)}}{f}}-\beta_{i,j}
            \right)^2
        +\beta_0.
    \end{align*}
    
    By construction, we have $(\PhiCphi(S_Q))_\Aff= X_\Aff \cap Q_\beta$. We may write $d_\beta$ as a quotient of polynomials
 \[d_\beta=\frac{g_\beta}{D}\quad
 \text{ with }\quad D:={\prod_{i=1}^n \left(\dotPoduct{C_i^{(0)}}{f}\right)^2}
 \]
    where 
    \[g_\beta= 
    \beta_0 D+\sum_{i=1}^n\frac{D}{\left(\dotPoduct{C_i^{(0)}}{f}\right)^2 }\sum_{j=1}^h 
        \left({(\dotPoduct{C_i^{(1)}}{f})}-\beta_{i,j}{(\dotPoduct{C_i^{(0)}}{f}})\right)^2
    \]
    is a homogeneous degree $2\ee n$ polynomial in $s,t$. Because $\beta_0$ is generic, 
    if $g_\beta$ has no common roots with $D$ then $g_\beta$ has $2\ee n$ distinct roots. 
    On the other hand, if $g_\beta$ had a common root with $D$, then there would exist a unique $i=a$ such that $\left(\dotPoduct{C_i^{(0)}}{f}\right)$ 
    vanishes on that root.
    This implies the polynomial
    \[
            \sum_{j=1}^\dimPPImage 
        {\left(\dotPoduct{C_a^{(j)}}{f}\right)}^2
    \]
    vanishes at the root. By the genericity of $\dotPoduct{C_a^{(j)}}{f}$ this cannot happen, thus $g_\beta$ has $2\ee n$ roots. 
      
    The degree of $g_\beta$ equals $\chi(X\cap Q_\beta).$ 
    It can be obtained by counting the degree of, for instance, the first term as we clear denominators. Multiplying through by $\dotPoduct{C_1^{(0)}}{f}^2$ gives a leading monomial of degree $2\ee n$, to which multiplication with subsequent $\dotPoduct{C_i^{(0)}}{f}^2$ adds $2\ee n$ for each pair of terms from the same camera. Thus, in an $n$-camera configuration, $\chi(X_\Aff\cap Q_\beta)=2\ee n$. 

  If $X_\Aff$ is smooth,
    an application of \cite[Theorem 3.8]{MRW-multiview}
    implies that 

    \begin{equation}\label{eq:smooth-curve-topological-formula}
        \affEDdeg{X}=-(\chi(X)-\chi(X\cap\HInfinity)-\chi(X\cap Q_\beta)).      
    \end{equation}
    Hence
    $\affEDdeg{X}=-(\chi(\PhiCphi(\PP^1))-
            \chi(\PhiCphi(S_\infty))-\chi(\PhiCphi(S_Q)))$,   
    and therefore
    \begin{equation}
        \affEDdeg{X} = (-1)(2-\ee n-2 \ee n)= 3\ee n-2.
    \end{equation}
\smallskip

If $X_\Aff$ 
 is singular then
we use \Cref{theorem:ed-affine}, which implies 
for general $\boldsymbol{\beta}:=(\beta_0,\beta_1,\dots,\beta_{\dimPPImage n})\in \mathbb{C}^{\dimPPImage n+1}$ that
    $\affEDdeg{X}$ equals
    \[ 
        -\left(
        \chi(X_{\reg})-\chi(X_{\reg}\cap H_\infty)-\chi((X_\reg)_\Aff\cap Q_\beta) 
        +\sum_{x\in \Sing(X_\Aff)\setminus Q_\beta} m_x
        \right)
        \]
    where 
    $X_\reg$ denotes the smooth points of $X$ and
    $m_x$ denotes the multiplicity of $x$.
        Similar arguments as in the smooth case show $\chi(X_\reg\cap \HInfinity)=\ee n$ and $$\chi((X_\reg)_\Aff\cap Q_\beta)=2\ee n$$
       as a consequence of the genericity of $\camConf$.
       Also by the genericity of $\camConf$ 
       we have
       \[
       \sum_{x\in \Sing(X_\Aff)\setminus Q_\beta} m_x=\sum_{x\in \Sing(X)} m_x.\]
Substituting in the known quantities, we find
 \begin{equation}\label{eq:need-singularity-cancellation}
    \affEDdeg{X}= 
        (-1)\left(
        \chi(X_{\reg})-\ee n-2\ee n 
        +\sum_{x\in \Sing(X)} m_x
        \right).
 \end{equation}       
Since each point in the singular locus of $X$ is a node and $X=\PhiCphi(\PP^1)$, we have 
\[
    \chi(X_\reg)=    \chi(\PP^1)
    -
    \sum_{x\in \Sing(X)} m_x.
    \]
Plugging into \eqref{eq:need-singularity-cancellation}, we get the desired formula. 

By \Cref{lemma:smooth-or-nodal}, the result follows.
\end{proof}

    \begin{corollary}\label{corollary:n-2-is-enough}    
    Suppose $\parametricMap$ is as in \eqref{eq:phi-P1-to-Pell} with $N\geq3$ and $\dimPPImage\geq2$.  Let $\irreducibleCameraSpace\subset\PP^{(\dimPPImage+1)\times (\dimPPWorld+1)-1}$ be an irreducible projective variety. 
    If $C_1,\dots,C_n$ are generic points in $\irreducibleCameraSpace$ and
    \begin{equation}\label{eq:ed-degree-C-family}
    \affEDdeg{\snap{(C_1,\dots,C_n)}{\phiWorldVariety}}=3\ee n-2 
    \end{equation}
    holds for $n=1,2$, 
    then \eqref{eq:ed-degree-C-family}
    holds for all $n\geq 1$.
    \end{corollary}

\begin{proof}
    When $\snap{(C_1,\dots,C_n)}{\phiWorldVariety}$ is smooth this follows from the proofs of \Cref{lemma:smooth-or-nodal} and \Cref{thm:ed-curve}, as each genericity condition only involves a single camera or a pair of cameras. 
\end{proof}

\begin{example}[A non-generic family]

 Let $\worldVariety$ be the rational curve in $\PP^{5}$ of degree $5$
    given by 
    \[
        [s:t]\overset{\parametricMap}{\mapsto} [t^5:st^{5-1}:\dots:s^5].
        \]
    Define the following irreducible variety of cameras:
    \[\irreducibleCameraSpace:=\left\{
        \begin{bmatrix}
            c_1&c_2&0&0&0&0\\
            0&0&c_3&c_4&0&0\\
            0&0&0&0&c_5&c_6
        \end{bmatrix}\in \field^{3\times 6}: c_i\in\field
        \right\}.
        \]
If $C$ is a generic element in $\irreducibleCameraSpace$,
then  $\affEDdeg{\snap{C_1}{Y}}=9$.
{We see that this disagrees with the numbers in \eqref{eq:ed-degree-C-family}}.

\end{example}

\begin{example}[Non-generic two views]
    
Consider the following family of $3\times4$ camera matrices:
    \[\irreducibleCameraSpace=\left\{
        [c_{i,j}]\in \jose{\PP^{3\times 4-1}}: c_{1,1}=0
        \right\}.
        \]
Let $C_1$ and $C_2$ be generic points of $\irreducibleCameraSpace$.
If $\worldVariety$ is the rational normal curve in $\PP^3$ parameterized by $[t^3:st^2:s^2t:s^3]$,
then 
\[\affEDdeg{\snap{C_1}{\worldVariety}}=7
\text{ and }
\affEDdeg{\snap{(C_1,C_2)}{\worldVariety}}=13.\]
The single view case agrees with the ED degrees given by the formula in \Cref{thm:ed-curve}, but the 2-view case does not. 
\end{example}

\begin{example}[Smooth to singular]
   
    \newcommand{\efour}{4}
    Let $\worldVariety$ be the rational curve in $\PP^{\efour}$ of degree $\efour$
    given by 
    \[
        [s:t]\overset{\parametricMap}{\mapsto} [t^\efour:st^{\efour-1}:\dots:s^\efour].
        \]
    Let $C$ be a general $3\times {(\efour+1)}$ camera matrix.
    Then $\snap{C}{X}$ is a degree $\efour$ 
    curve in $\PP^2$ with 
    $\binom{\efour-1}{2}$ nodal singularities.
    The set of smooth points of $\snap{C}{X}\subset \PP^2$ has Euler characteristic  $\chi(\PP^1\setminus \{2\cdot \binom{\efour-1}{2} \text{ points}\})=2-2\cdot\binom{\efour-1}{2}$.
    The ED degree of $(\snap{C}{X})_{\operatorname{Aff}}$ is $3\cdot\efour-2$.    
    The projective ED degrees of $X$ and $\snap{C}{X}$ 
    both equal {$3\cdot\efour-2$}. As affine cones in $\CC^{\efour+1}$ and $\CC^3$, these varieties are parameterized by degree $\efour$ polynomials and one might predict from \cite[Section 2]{DHOST2016-ed-degree} that the ED degree is $(2\cdot \efour-1)^2$. The reason we do not get this bound is because these are general homogeneous polynomials of degree $\efour$ and not general degree $\efour$ polynomials.  
\end{example}

\begin{example}[Cuspidal cubic]
If $\parametricMap(s,t)=[s^3:st^2:t^3]$
and the camera configuration $\camConf$ consists of a single generic camera, 
then the ED degree of $(\snap{\camConf}{\phiWorldVariety})_\Aff$ is six.   
The cubic curve $\phiWorldVariety$ has a cusp as a singularity. Thus, $\chi(X_\reg)=\chi(\PP^1)-1$ and this is the reason for the drop in the ED degree from the ED degree of a cubic curve in $\PP^2$ with a node. 
\end{example}

\begin{remark}[Intersection theory and multidegrees]\label{remark:multidegree-perspective}
    The results in this section have an intersection theoretic~\cite{Fulton-intersection-theory} interpretation through the multidegrees presented in \Cref{ss:multiprojective}.
    Namely, the proof in \Cref{thm:ed-curve} is showing that the curve $\snap{\camConf}{\phiWorldVariety}\subset(\oneImageAmbient)^n$ has multidegree 
    $\sum_{i=1}^n \ee \cdot {T_1^{\dimPPImage-1}T_2^{\dimPPImage-1}\cdots T_n^{\dimPPImage-1}}$ and intersects transversally each of the  hyperplanes~$V(x_{i,0})$ for $i=1,\dots,n$.  This implies $\chi(X\cap H_\infty)=\ee n$.
    
    In addition, the proof shows that
    the hypersurface defined by 
    \[q_\beta:= \sum_{i=1}^n \sum_{j=1}^\dimPPImage(x_{i,j}-\beta_{i,j}x_{i,0})^2+\beta_0\prod_{i=1}^nx_{i,0}^2\]    
    (obtained by homogenizing the polynomial in \eqref{eq:Q-beta-U-chart})
    also transversally intersects the curve $\snap{\camConf}{\phiWorldVariety}$ with no solutions at infinity.
    Since the multidegree of $V(q_\beta)$ is 
    $\sum_{i=1}^n 2T_i$, we have  $\chi(X\cap V(q_\beta))=2\ee n$ because 
    \[\left(\sum_{i=1}^n 2T_i\right)\cdot \left(\sum_{i=1}^n \ee \cdot {T_1^{\dimPPImage-1}\cdots T_n^{\dimPPImage-1}}\right)= 2\ee n\cdot T_1^\dimPPImage\cdots T_n^\dimPPImage \in \mathbb{Z}[T_1,\dots,T_n]/\langle T_1^{\dimPPImage+1},\dots, T_n^{\dimPPImage+1}\rangle.\]

\end{remark}

\section{Settling conjectures of Duff--Rydell}\label{s:ed-degree-conjectures}

In this section, we present our main application of our \jose{\Cref{thm:ed-curve}}.
\Cref{ss:anchored-world-schubert} defines our main varieties of interest: a {set of} multiview varieties anchored at a subvariety of the Grassmannian. 
\Cref{ss:exterior-wedge-cameras} gives a complementary approach for obtaining these multiview varieties using wedge cameras. It also provides the necessary exterior algebra background. 
Using our results in \Cref{s:multiview-curves}, we prove a pair of conjectures 
in \Cref{ss:equal-ED-degrees}.

\subsection{Anchoring the $n$-view variety 
}\label{ss:anchored-world-schubert}

{\Cref{def:MultiviewVariety-V2} considers  
 the {point multiview variety} 
$\snap{\camConf}{\worldVariety}$
 anchored at $\worldVariety$
 as the Zariski closure of the image of the map \eqref{eq:PhiCam-V2} for $\worldVariety\subset\worldAmbient$.
  }
In this section, we carry this idea over to a  map between Grassmannians. 
The associated multiview variety (called a \emph{generalized multiview variety} in \cite{duff-rydel-2024metricmultiviewgeometry}) is obtained as the Zariski closure of a rational map
\begin{equation}\label{eq:gr-to-gr-product-camera}
\begin{aligned}
    \Phi_{\camConf,k}:\Gr(k,\worldAmbient) &\dashrightarrow \Gr(k,\oneImageAmbient)^n\\
    P &\mapsto (C_1 \cdot P, \dots, C_n \cdot P)
\end{aligned}
\end{equation}
for a camera arrangement $\camConf=(C_1,\dots,C_n)$ of full-rank $(\dimPPImage+1)\times (\dimPPWorld+1)$ matrices.
Here, the multiplication $C_i\cdot P$ denotes the $k$-plane ($k$-dimensional linear subspace of $\worldAmbient$) defined as the span of $C_iX_0,\dots,C_iX_k$ for $P$ spanned by $X_0,\dots,X_k$. 
If $k=0$, then $\Phi_{\camConf,k}$ agrees with \eqref{eq:PhiCam-V2}.

When the domain of $\Phi_{\camConf,k}$ 
is restricted to a subvariety, the Zariski closure of the image is {still} referred to as an {anchored multiview variety}.
Given a 
subvariety
$\Lambda\subseteq\Gr(k,\worldAmbient)$ of the Grassmannian, one obtains a multiview variety as the Zariski closure of the image of 
\begin{equation}\label{cam-lambda}
\begin{aligned}
    \Lambda &\dashrightarrow \Gr(k,\oneImageAmbient)^n\\
    P &\mapsto (C_1 \cdot P, \dots, C_n \cdot P).
\end{aligned}
\end{equation}

We summarize this discussion with the following definition.

\begin{definition}\label{definition:schubert-multiview}
Let $\Lambda$ be a subvariety of $\Gr(k,\worldAmbient)$.
For a camera arrangement $\camConf=(C_1,\dots,C_n)$ of full-rank $(\dimPPImage+1)\times (\dimPPWorld+1)$ matrices,
we denote the Zariski closure of the image of \eqref{cam-lambda}
by
$\snap{\camConf}{\worldSchubert}$ and call this the \demph{$n$-view variety anchored at $\worldSchubert$} with respect to $\camConf$.
When $k=1$, we say it is a \demph{line $n$-view variety}.
\end{definition}

\begin{remark}\label{rmk:duff-rydell-notation}
    In \cite{duff-rydel-2024metricmultiviewgeometry},
    the full data of a multiview variety is recorded in the notation $\mathcal{M}_{\camConf,k}^{\Lambda,\dimPPWorld,\dimPPImage}$. We opt for the more parsimonious {notation}
    and refer to the catalogue in \cite[Theorem 3.3]{duff-rydel-2024metricmultiviewgeometry} for a complete description of multiview varieties arising from projections of points and lines in 1, 2, and 3-dimensional
    projective space.
\end{remark}

As defined, $\snap{\camConf}{\Lambda}$ is a subvariety of a product of Grassmannians. 
Next, we embed $\snap{\camConf}{\Lambda}$ into a product of projective spaces. 
This is done using the  Pl\"ucker embedding  on  a product of Grassmannians,
\begin{equation}\label{eq:khn-plucker-Gr-k-h}
\demph{{$\plucker_{k,h,n}$}}
:(\Gr(k,\oneImageAmbient))^n\to \left(\PP^{\binom{\dimPPImage+1}{k+1}-1}\right)^n.
\end{equation}
{When the context is clear, we suppress the dependence of the map $\iota$ on $k$, $\dimPPImage$, and~$n$.}
Thus {$\plucker(\snap{\camConf}{\Lambda})_{\operatorname{Aff}}$} 
is an affine variety contained in 
$\left(\CC^{\binom{\dimPPImage+1}{k+1}-1}\right)^n$.

\begin{example}[Defining $\Xhn{\dimPPImage}{n}$]\label{example:L3-line-view}
    Recall $\LThree\subset\Gr(1,\PP^3)$ from \Cref{ss:schubert}.
    Suppose $\dimPPImage=2$ or $\dimPPImage=3$.
    If $\hCamConf$ is an arrangement of $n$-cameras of size $(\dimPPImage+1)\times 4$,
    then $\snap{\hCamConf}{\LThree}$ is a subvariety of $\Gr(1,\oneImageAmbient)^n$.
    Taking the image of the Pl\"ucker embedding
    \[\plucker_{1,h,n}:\snap{\hCamConf}{\LThree}\to (\PP^{\binom{h+1}{2}-1})^n
    \]
    we obtain a multiprojective variety, which we denote by \demph{$\Xhn{\dimPPImage}{n}$}. 
    \Cref{table:small-cases} highlights properties of 
    {$\Xhn{\dimPPImage}{n}$}
    for small $n$.
\end{example}

     \begin{table}[htb!]
         \centering
         \begin{tabular}{c|c|c|c}
            & Ambient space  
            & Multidegree 
            & $\affEDdeg{\Xhn{\dimPPImage}{n}}$ \\
            \hline
    $\dimPPImage=2,\,n=1$&  
        $\PP^{2}$   
        & $2T_1$
        & 4
        \\ 
    $\dimPPImage=2,\,n=2$
    &  $\PP^{2}\times \PP^{2}$
        & $\jose{(2T_1+2T_2)T_1T_2}$
        & 10
        \\
    $\dimPPImage=3,\,n=1$& 
        $\PP^{5}$   
        & $\jose{2T_1^4}$
        & 4
        \\ 
    $\dimPPImage=3,\,n=2$
        &  $\PP^{5}\times \PP^{5}$
        & $\jose{(2T_1+2T_2)T_1^4T_2^4}$
        & 10
    \end{tabular}
         \caption{Summary of $\snap{\hCamConf}{\LThree}$ with $\hCamConf$ having $n=1,2$ cameras in its arrangement.
         }
         \label{table:small-cases}
     \end{table}
    
\medskip

\subsection{Exterior algebra review and wedge cameras}\label{ss:exterior-wedge-cameras}

Next, we show how to realize
$\snap{\hCamConf}{\LThree}$ as a point multiview variety anchored at $\iota(\LThree)\subset \PP^5$.
To that end, we review the exterior algebra {(see \cite[Chapter XVI]{MLB1988-algebra} for a textbook reference)} to state
\Cref{definition:wedge-cameras}.

\newcommand{\kWedgeMap}[2]{f_{\wedge^{#1}#2}}

\newcommand{\kWedgeCamera}[2]{{\wedge^{#1}#2}}

\newcommand{\NOne}{N}
\newcommand{\NTwo}{M}

Recall the \emph{exterior algebra} of a vector space $V$ over a field $\mathbb{K}$ is the graded algebra
\[
\bigwedge V := \bigoplus_{k=0}^{\dim V} \bigwedge^{k} V,
\]
where $\bigwedge^k V$ is the space of \emph{$k$-vectors}, the dual of the space of alternating $k$-linear forms on $V$. 
The wedge product satisfies
$u \wedge v = - v \wedge u$ and  $v \wedge v = 0$.
A  \emph{simple wedge} is one of the form
$
v_1 \wedge \cdots \wedge v_k$, $v_i \in V.
$
Every element of $\bigwedge^k V$ is a linear combination of simple wedges. 
The dimension of $\bigwedge^k V$ is $\binom{\dim V}{k}$, and a basis is given by
\[
e_{i_1} \wedge e_{i_2} \wedge \cdots \wedge e_{i_k},\quad i_1 < i_2 < \cdots < i_k
\]
where $\{e_i\}$ is a basis for $V$.
Any linear map $f: V \to W$ induces the map
\begin{equation}\label{eq:linear-map-wedge}
\bigwedge^k f : \, \bigwedge^k V \longrightarrow \bigwedge^k W,\qquad
\bigwedge^k f(v_1 \wedge \cdots \wedge v_k) = f(v_1) \wedge \cdots \wedge f(v_k).
\end{equation}

\smallskip
In our vision setting, each camera matrix $C$ is a linear map
$f_C: \mathbb{K}^{\dimPPWorld+1} \to \mathbb{K}^{\dimPPImage+1}.$
After fixing ordered bases of $\field^{\dimPPWorld+1}$ and $\field^{\dimPPImage+1}$, we obtain induced wedge bases on
$\bigwedge^{k}\field^{\dimPPWorld+1}$ and
$\bigwedge^{k}\field^{\dimPPImage+1}$. 
The \demph{wedge camera matrix $\kWedgeCamera{k}{C}$} represents the linear map 
\[
\bigwedge^k f_C:\;
\bigwedge^{k}\field^{\dimPPWorld+1}\to
\bigwedge^{k}\field^{\dimPPImage+1},
\text{ where }
v_1 \wedge \cdots \wedge v_k
\mapsto f_C(v_1) \wedge \cdots \wedge f_C(v_k).
\]
The entries of $\kWedgeCamera{k}{C}$ consist of all $k\times k$ minors of $C$, arranged according to the chosen multi-index ordering.

\begin{definition}\label{definition:wedge-cameras}
 
If $\camConf=(C_1,\dots,C_n)$ is an arrangement of $n$ cameras, then we call
$
\left(\kWedgeCamera{k+1}{C_1},\dots,\kWedgeCamera{k+1}{C_n}\right)
$
 the \demph{$(k+1)$-wedge camera arrangement of $\camConf$}.
\end{definition}

\begin{example}\label{example:wedge-camera-3-6}

Suppose $C$ is the $3\times4$ matrix
\[
C = \begin{bmatrix}
1 & 2 & 3 & 4\\
0 & -1 & 0 & 0\\
0 & 0 & 5 & 0
\end{bmatrix}.
\]
Then 
$\kWedgeMap{2}{C}:\; \bigwedge^2\mathbb{K}^4\to\bigwedge^2\field^3$
has matrix
$\kWedgeCamera{2}{C}\in\field^{3\times 6}$ with entries being the $2\times 2$ minors of $C$.

 We use the reverse lexicographic ordering 
$((1,2),(1,3),(2,3))$ on rows and \newline
 $((1,2),(1,3),(2,3),(1,4),(2,4),(3,4))$ for columns.
Then, 
\[
\kWedgeCamera{2}{C} =
\begin{bmatrix}
-1 & 0 & 3 & 0 & 4 & 0\\
0 & 5 & 10 & 0 & 0 & -20\\
0 & 0 & -5 & 0 & 0 & 0
\end{bmatrix}.
\]
Each entry is $\det C_{I,J}$ with $I$ a row pair and $J$ a column pair in the stated order.
\end{example}

The proof of the following proposition is sketched in \cite{duff-rydel-2024metricmultiviewgeometry}. It is a consequence of standard results from multilinear algebra. 

\begin{proposition}\label{prop:revised-grassmannian-multiview}
Let $\dimPPImage\geq k$.
Fix a camera arrangement $\camConf=(C_1,\dots,C_n)$ of full-rank $(\dimPPImage+1)\times (\dimPPWorld+1)$ matrices and a subvariety $\Lambda$ of $\Gr(k,\worldAmbient)$.
Then 
\[\grCamConf=(\wedge^{k+1}C_1, \dots, \wedge^{k+1}C_n)\] 
is a {wedge} camera arrangement of $\binom{\dimPPImage+1}{k+1}\times\binom{\dimPPWorld+1}{k+1}$ full rank matrices. 
Moreover, 
$\plucker_{k,h,n}(\snap{\camConf}{\Lambda})$
is a point multiview variety anchored at $\plucker_{k,N,1}(\Lambda)$: 
\[\plucker_{k,h,n}(\snap{\camConf}{\Lambda})=\snap{\grCamConf}{\plucker_{k,N,1}(\Lambda)}.
\]
\end{proposition}

As a consequence of \Cref{prop:revised-grassmannian-multiview} we get the commutative diagram in \Cref{fig:all-the-maps-schubert}. 

\begin{figure}[htb!]
     \centering
\tikzcdset{
  every cell/.append style = {font = \large},   
  every label/.append style = {font = },  
}
\tikzset{
  >=Latex,                     
  cdArrow/.style = {line width=1.5pt}, 
}
\begin{tikzcd}[
  column sep=large,         
  row sep=huge,            
  cells={nodes={inner sep=8pt}}
]
  \worldSchubert\subset\Gr(k,\worldAmbient) \arrow[r, "\Phi_{\camConf}"] \arrow[d, "\plucker_{k,\dimPPWorld,1}"'] 
  & \snap{\camConf}{\worldSchubert} \arrow[d, "\plucker_{k,\dimPPImage,n}",shift right=0pt] \arrow[r, ""]&\Gr(k,\oneImageAmbient)^n \arrow[d,"\plucker_{k,\dimPPImage,n}"]
\\
  \plucker_{k,\dimPPWorld,1}(\Lambda)\subset \PP^{\binom{\dimPPWorld}{k}} \arrow[r, "\Phi_{\grCamConf}"'] 
  &\snap{\grCamConf}{\plucker_{k,\dimPPImage,n}(\worldSchubert)}\arrow[r, ""]&(\PP^{\binom{\dimPPImage+1}{k}-1})^n
\end{tikzcd}
    \caption{This commutative diagram shows two approaches for obtaining the $k$-plane multiview variety anchored at $\Lambda$.
    }
    \label{fig:all-the-maps-schubert}
\end{figure}

\subsection{Resolving a pair of
ED degree conjectures}\label{ss:equal-ED-degrees}

Recall the subvariety $\LThree$ of $\Gr(1,\PP^3)$ 
from \Cref{ss:schubert}
and  \Cref{definition:wedge-cameras}.
We focus on $\snap{\hCamConf}{\LThree}$,
the line view variety anchored at $\LThree$ with respect to a generic arrangement $\hCamConf$ of $n$ cameras of size $(\dimPPImage+1)\times4$.
Throughout this subsection,
we denote the Pl\"ucker embedding of  
$\snap{\hCamConf}{\LThree}$
by
\[\Xhn{\dimPPImage}{n}:=\plucker_{1,\dimPPImage,n}(\snap{\hCamConf}{\LThree}).\] 

\begin{example}\label{example:L3-point-view}
Continuing with \Cref{example:L3-line-view}, 
suppose $\hGCamConf$
is the $2$-wedge camera arrangement of~$\hCamConf$.
Then $\hGCamConf$ is an arrangement of cameras of size $\binom{h+1}{2}\times 6$. 
By \Cref{prop:revised-grassmannian-multiview},
we have
\[  
    \Xhn{\dimPPImage}{n}
    =
    \snap{\hGCamConf}{\plucker_{1,3,1}(\LThree)}.
    \]
Furthermore, 
combining the parameterization in \Cref{ex:embedding-L3} with the results of 
\Cref{prop:revised-grassmannian-multiview} yields the commutative diagram in \Cref{fig:all-the-maps-L-three}.
The commutative diagram 
shows  $\Phi_{\hGCamConf}\circ f$ is a parameterization of
the line $n$-view variety anchored at~$\LThree$.
\end{example}

\begin{figure}[hbtb]
    \centering
\tikzcdset{
  every cell/.append style = {font = \large},   
  every label/.append style = {font = },  
}
\tikzset{
  >=Latex,                     
  cdArrow/.style = {line width=1.5pt}, 
}
\begin{tikzcd}[
  column sep=large,         
  row sep=huge,            
  cells={nodes={inner sep=8pt}}
]
  & \LThree\subset\Gr(1,\PP^3) \arrow[r, "\Phi_{\hCamConf}"] \arrow[d, "\plucker_{1,3,1}"'] 
  & \snap{\hCamConf}{\LThree} \arrow[d, "\plucker_{1,\dimPPImage,n}",shift right=0pt] \arrow[r, ""]&\Gr(1,\oneImageAmbient)^n \arrow[d,"\plucker_{1,h,n}"]
\\
\PP^1 \arrow[r, "f"'] 
  & \plucker_{1,3,1}(\LThree)\subset \PP^5 \arrow[r, "\Phi_{\hGCamConf}"'] 
  &\snap{\hGCamConf}{\plucker_{1,3,1}(\LThree)}\arrow[r, ""]&(\PP^{\binom{\dimPPImage+1}{2}-1})^n
\end{tikzcd}
    \caption{For $h=2,3$, this shows the equality 
    $\plucker_{1,\dimPPImage,n}(\snap{\hCamConf}{\LThree})=\snap{\hGCamConf}{\plucker_{1,3,1}(\LThree)}$.
    The unlabeled horizontal arrows are inclusions and $\plucker_{k,h,n}$  is as defined in \eqref{eq:khn-plucker-Gr-k-h}.}
    \label{fig:all-the-maps-L-three}
\end{figure}

\begin{theorem}[\cite{duff-rydel-2024metricmultiviewgeometry} Conjecture 7.4.5 and Conjecture 7.4.6]\label{thm:duff-rydell-conj-V2}
    For $\dimPPImage=2,3$,
    let $\hCamConf$ be a generic configuration of $n$ cameras of size $(\dimPPImage+1)\times4$.    %
    If $\Xhn{n}{h}$ denotes $\plucker(\snap{\hCamConf}{\LThree})$,
    then
\begin{align}\label{eq:6n-2}
    \affEDdeg{\Xhn{\dimPPImage}{n}}&=6n-2.
    \end{align}
\end{theorem}    

\begin{proof}
\jose{Let $\dimPPImage=2$ or $\dimPPImage=3$, and}
    let $\hGCamConf=(D_1,\dots,D_n)$ denote the $2$-wedge camera configuration of $\hCamConf$.

    By \Cref{prop:revised-grassmannian-multiview} and \Cref{example:L3-point-view}, 
    we have
    \begin{align*}
        \plucker_{1,\dimPPImage,n}(\snap{\hCamConf}{\LThree})
        \;&=\; \snap{\hGCamConf}{\plucker_{1,3,1}(\LThree)}\\
        \Xhn{\dimPPImage}{n}
        \;&=\; \snap{\hGCamConf}{\plucker_{1,3,1}(\LThree)}.
    \end{align*} 
    Since $\hCamConf$ is generic, each $D_i$ is a generic point in the irreducible variety of wedge cameras of size $\binom{h+1}{2}\times 6$.
    {It follows from \Cref{fig:all-the-maps-L-three}
    that $\plucker_{1,3,1}(\LThree)$ is parameterized by $f$ in \Cref{ss:schubert}.}
    Thus, by \Cref{corollary:n-2-is-enough},
    it suffices to show \eqref{eq:6n-2} for $n=1,2$.
    This holds by \Cref{table:small-cases} and the result follows.
\end{proof}

\newcommand{\iPointOne}[1]{p_{1,#1}}
\newcommand{\iPointTwo}[1]{p_{2,#1}}
\newcommand{\bezierWorld}[1]{\Lambda_{#1}}

\section{One-parameter families of 3D lines}\label{sec:One-parameterFamilies}

We study an additional model to further motivate the investigation of ED Degrees of subvarieties of the Grassmannian.
One important example comes from rational scrolls. These are surfaces that represent one-parameter families of 3D lines and arise naturally in computer vision applications. In this section, we consider an instance of such a family: lines sweeping along B\'ezier curves to provide a one-parameter line-view variety. This is obtained by mapping a subvariety of $\Gr(1,\PP^3)$ to $\Gr(1,\PP^2)^n$.

   \medskip
For background \jose{(e.g., \cite[Chapter 5]{Goldman-pyramid-algorithms})}, recall a \demph{B\'ezier  curve}
in $\field^\dimPPWorld$ of degree 
$\ee$ is defined by a set of $\ee+1$ control points~
\[P_0,\dots,P_\ee \in\field^{\dimPPWorld}\] with $N\geq 2$.
    The curve is given parametrically by
\begin{equation}\label{eq:Bezier-affine}
\mathbf{B}(1,t) = \sum_{i=0}^{\ee} 
 B_{i,\ee}(1,t)\cdot {P}_i, \quad  \, t\in[0,1]
\end{equation}
where $B_{i,\ee}(s,t)$ is the $i$th homogenized Bernstein polynomial of degree $\ee$,
\[B_{i,\ee}(s,t) := \binom{\ee}{i} (s - t)^{\ee - i} t^i.\]
This parameterization  has the advantage that the curve segment 
$\{\mathbf{B}(1,t) :  t\in[0,1]\}$ on the standard affine chart in $\field^{\dimPPWorld}$ is in the convex hull of the control points, starts at $P_0=\mathbf{B}(1,0)$, and ends at $P_\ee=\mathbf{B}(1,1)$. 
Using the homogenized polynomials, the parameterization from \eqref{eq:Bezier-affine} extends to a  
projective map
\[\PP^1 \longrightarrow \PP^{\dimPPWorld},\qquad 
[s:t]\longmapsto [s^\ee:\mathbf{B}(s,t)].\]

\newcommand{\fBB}{f_{\mathbf{B}_1,\mathbf{B}_2}}

\begin{figure}[hbt!]
    \centering
        \includegraphics[width=0.42\linewidth]{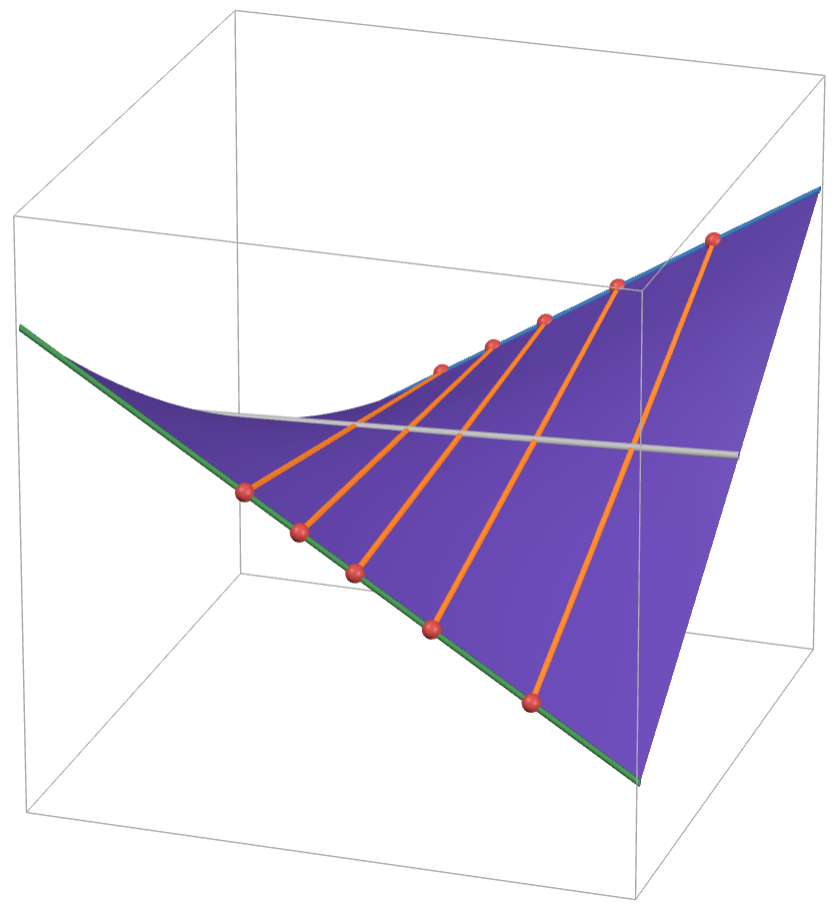} \,\,\,\,\,\,\,\,\,
    \includegraphics[width=0.5\linewidth]{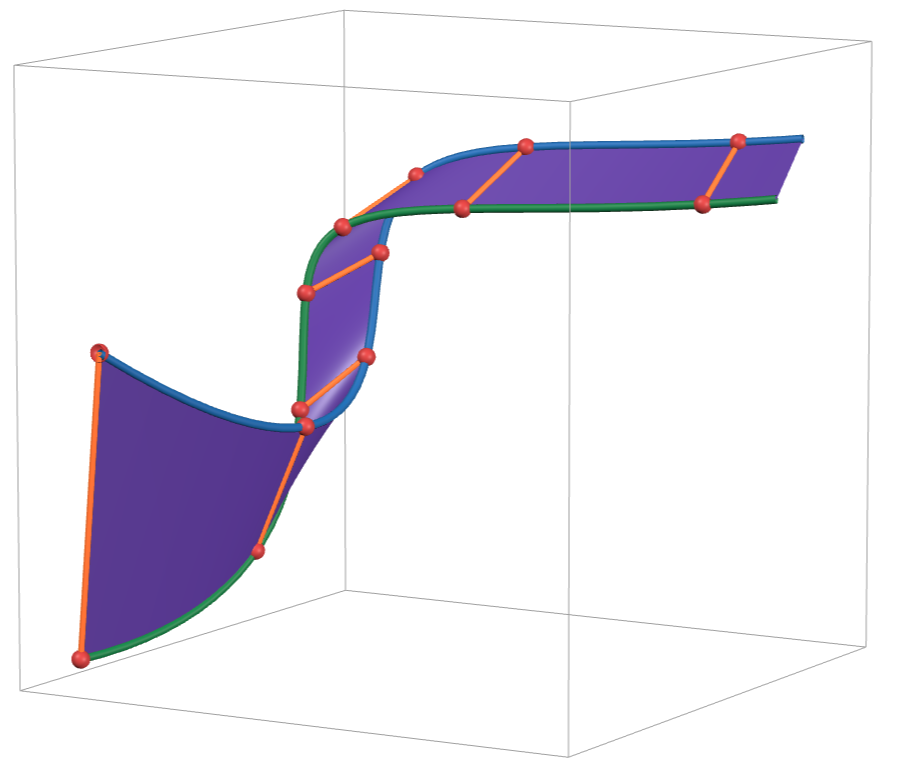}
    \caption{The purple ruled surfaces consist of line segments connecting a pair of B\'ezier curves of degree one (left) and four (right). The lines spanned by the segments in a ruling give
    $\bezierWorld{\mathbf{B}_1,\mathbf{B}_2}
    \subset\Gr(1,\PP^3)$. The left surface has two rulings as discussed in~\Cref{ss:schubert}.
    }
    \label{fig:Bezier-scrolls}
\end{figure}

\smallskip
We now construct a one-parameter family of lines.
For simplicity, we assume $\dimPPWorld=3$.
Suppose 
we have a pair of B\'ezier curves 
 $\mathbf{B}_1(1,t)$
 and 
 $\mathbf{B}_2(1,t)$ 
 in $\field^3$ with degree $\ee_1$ and $\ee_2$, respectively. 
Their projective closures are the rational curves 
\[
X_i:=\left\{
    [s^{\ee_i}: \mathbf{B}_i(s,t)]\in \PP^3: [s:t]\in \PP^1\right\}.
\]
Under genericity assumptions of the control points, we may assume $X_1\cap X_2$ is empty. Thus, the points
\[
[s^{\ee_1}: \mathbf{B}_1(s,t)] \text{ and } [s^{\ee_2}: \mathbf{B}_2(s,t)]
\]
are distinct 
for all $[s:t]\in \PP^1$, 
and span a line in $\PP^3$. 
This set of lines {yields} a subvariety of $\Gr(1,\PP^3)$,
{which we denote by} $\bezierWorld{\mathbf{B}_1,\mathbf{B}_2}$.
This is visualized in \Cref{fig:Bezier-scrolls}.

Let $\irreducibleCameraSpace\subset\PP^{11}$ be an irreducible projective variety of $3\times 4$ camera matrices. 
If $C_1,\dots,C_n$ are generic points in $\irreducibleCameraSpace$, our line-multiview variety of interest is 
\[\snap{(C_1,\dots,C_n)}{\bezierWorld{\mathbf{B}_1,\mathbf{B}_2}}\;\subset\; \Gr(1,\PP^2)^n.\]
Recalling \Cref{definition:schubert-multiview}, we say it is anchored at $\bezierWorld{\mathbf{B}_1,\mathbf{B}_2}$ with respect to the camera arrangement $\camConf=(C_1,\dots,C_n)$.
\medskip

\newcommand{\XBBn}{X_{\mathbf{B}_1,\mathbf{B}_2,n}}

\begin{theorem}\label{theorem:bezier}    
    Let $\irreducibleCameraSpace\subset\PP^{3\times 4-1}$ be an irreducible projective variety, and fix generic points
     $C_1,\dots,C_n$ in $\irreducibleCameraSpace$.
     For B\'ezier curves 
     $\mathbf{B}_1$ and  $\mathbf{B}_2$ with $\ee_1+1$ and $\ee_2+1$ generic control points,
     let $\XBBn$ denote the multiprojective variety 
     \[\plucker\left(
        \,\snap{(C_1,\dots,C_n)}{\bezierWorld{\mathbf{B}_1,\mathbf{B}_2}}\,
     \right)\subset (\PP^2)^{n}.\] 
     If 
    \begin{equation}\label{eq:ed-degree-D-family}
    \affEDdeg{\XBBn}=3(\ee_1+\ee_2) n-2, 
    \end{equation}
    holds for $n=1,2$, 
    then \eqref{eq:ed-degree-D-family}
    holds for all $n\geq 1$.
\end{theorem}
\begin{proof}

Under the Pl\"ucker embedding a line in $\bezierWorld{\mathbf{B}_1,\mathbf{B}_2}$ is determined by the six maximal minors of the $2\times4$ matrix
\[
\begin{bmatrix}
s^{\ee_1}&\mathbf{B}_1(s,t)\\
s^{\ee_2}&\mathbf{B}_2(s,t)
\end{bmatrix}.
\]
This yields a degree $\ee_1+\ee_2$ rational map $\fBB:\PP^1\to \PP^5$, assuming the genericity of the control points of the B\'ezier curves.

By the commutative diagram in \Cref{ss:exterior-wedge-cameras}, which is a consequence of \Cref{prop:revised-grassmannian-multiview}, we have
\[
\XBBn=
{\snap{(D_1,\dots,D_n)}{\plucker(\bezierWorld{\mathbf{B}_1,\mathbf{B}_2})}}
\]
where $ D_i=\wedge^2C_i$.
\jose{Since ${\plucker(\bezierWorld{\mathbf{B}_1,\mathbf{B}_2})}=\fBB(\PP^1)$ is a rational curve in $\PP^3$ of degree $\ee_1+\ee_2$, the result follows from \Cref{corollary:n-2-is-enough} as the $D_i$ are generic points in the irreducible variety $\{\wedge^2C: C\in \irreducibleCameraSpace\}\subset\PP^{3\times 6-1}$.}
\end{proof}

\section{Conclusion and future directions}\label{s:conclusion}
In this article, we considered one-dimensional multiview varieties. We used topological arguments to give the first theoretical results for ED degrees of multiview varieties anchored at a Schubert variety. One of the key features of our proof is that each genericity condition is the complement of a union of hypersurfaces whose defining equations only rely 
\jose{on a single camera or a camera pair.}
This resulted in a surprising \Cref{corollary:n-2-is-enough} to determine the ED degrees of multiview varieties where cameras have a prescribed structure. This is relevant for families of ``dual cameras'' found in resectioning~\cite{MR4919571} as well as calibrated 
cameras~\cite{MR3705779,MR3636592} used in practice. 
{A natural extension is to generalize our results on the ED degree of curves to higher-dimensional varieties.}
\jose{We leave for future work the task of combining the techniques in \cite{MRW-multiview}, which were used to prove \eqref{eq:big-conjecture},
with the approaches we outlined in our proofs to give a generalization of our \Cref{corollary:n-2-is-enough} for all parameterized anchored mutliview varieties. 
 Another interesting problem is to determine the ED degree of line multiview varieties using alternative embeddings of the Grassmannian, such as those appearing in recent work on Grassmannian optimization~\cite{Lai-Lim-Ye-2025,LLY2025-gro,DFRS2025-two-lives}.} \bella{A natural next step is to choose a representation of the embedded multiview variety that yields good properties for the bundle adjustment problem. For example, we can ask how the study of representations for 3D lines in projective space in \cite[Section 3]{BartoliSturm} generalizes to other multiview varieties.} \jose{We hope that further developments along these lines will deepen the interplay between algebraic geometry, multilinear algebra, and computer vision.}

\subsubsection*{Acknowledgments}
We thank Timothy Duff for his comments on early versions of the~article. 

\subsubsection*{Data availability statement}The calculations supporting the results of this article can replicated using the open source software \texttt{Macaulay2}~\cite{M2} using the \texttt{EuclideanDistanceDegree.m2} package.
The code for these calculations is available at the~\href{https://github.com/JoseMath/EuclideanDistanceDegree/blob/multiview-tst-file/test-multiview.m2}{GitHub repository}.

\bibliographystyle{siam-no-dash-title-color-links} \bibliography{refs}

\begin{thebibliography}{10}

\bibitem{AST-multidegree-hilbert-scheme}
{\sc C.~Aholt, B.~Sturmfels, and R.~Thomas}, {\em \titleColor{A {H}ilbert scheme in computer vision}}, Canad. J. Math., 65 (2013), pp.~961--988.
\newblock \href{https://doi.org/10.4153/CJM-2012-023-2}{\texttt{DOI}}.

\bibitem{aluffi-harris-projective-ed}
{\sc P.~Aluffi and C.~Harris}, {\em \titleColor{The {E}uclidean distance degree of smooth complex projective varieties}}, Algebra Number Theory, 12 (2018), pp.~2005--2032.
\newblock \href{https://doi.org/10.2140/ant.2018.12.2005}{\texttt{DOI}}.

\bibitem{BartoliSturm}
{\sc A.~Bartoli and P.~Sturm}, {\em \titleColor{Structure-from-motion using lines: Representation, triangulation, and bundle adjustment}}, Computer Vision and Image Understanding, 100 (2005), pp.~416--441.
\newblock \href{https://doi.org/https://doi.org/10.1016/j.cviu.2005.06.001}{\texttt{DOI}}.

\bibitem{breiding2024line}
{\sc P.~Breiding, T.~Duff, L.~Gustafsson, F.~Rydell, and E.~Shehu}, {\em \titleColor{Line multiview ideals}}, Comm. Algebra, 52 (2024), pp.~4204--4225.
\newblock \href{https://doi.org/10.1080/00927872.2024.2343762}{\texttt{DOI}}.

\bibitem{BRST2023-line-view}
{\sc P.~Breiding, F.~Rydell, E.~Shehu, and A.~Torres}, {\em \titleColor{Line multiview varieties}}, SIAM J. Appl. Algebra Geom., 7 (2023), pp.~470--504.
\newblock \href{https://doi.org/10.1137/22M1482263}{\texttt{DOI}}.

\bibitem{MR4919571}
{\sc E.~Connelly, T.~Duff, and J.~Loucks-Tavitas}, {\em \titleColor{Algebra and geometry of camera resectioning}}, Math. Comp., 94 (2025), pp.~2613--2643.
\newblock \href{https://doi.org/10.1090/mcom/4030}{\texttt{DOI}}.

\bibitem{DFRS2025-two-lives}
{\sc K.~Devriendt, H.~Friedman, B.~Reinke, and B.~Sturmfels}, {\em \titleColor{The two lives of the {G}rassmannian}}, Acta Univ. Sapientiae Math., 17 (2025), pp.~Paper No. 8, 18.
\newblock \href{https://doi.org/10.1007/s44426-025-00007-x}{\texttt{DOI}}.

\bibitem{DHOST2016-ed-degree}
{\sc J.~Draisma, E.~Horobe\c~t, G.~Ottaviani, B.~Sturmfels, and R.~R. Thomas}, {\em \titleColor{The {E}uclidean distance degree of an algebraic variety}}, Found. Comput. Math., 16 (2016), pp.~99--149.
\newblock \href{https://doi.org/10.1007/s10208-014-9240-x}{\texttt{DOI}}.

\bibitem{duff-rydel-2024metricmultiviewgeometry}
{\sc T.~Duff and F.~Rydell}, {\em \titleColor{Metric multiview geometry -- a catalogue in low dimensions}}, 2024.
\newblock \href{https://arxiv.org/abs/2402.00648}{\texttt{arXiv:} \texttt{2402.00648}}.

\bibitem{eisenbud20163264}
{\sc D.~Eisenbud and J.~Harris}, \titleColor{{\em 3264 and All That: A Second Course in Algebraic Geometry}}, Cambridge University Press, 2016.
\newblock \href{https://books.google.com/books?id=9FZ6jwEACAAJ}{\texttt{URL}}.

\bibitem{EscobarKnutson2017-multidegree}
{\sc L.~Escobar and A.~Knutson}, \titleColor{{\em The Multidegree of the Multi-Image Variety}}, Springer New York, New York, NY, 2017, pp.~283--296.
\newblock \href{https://doi.org/10.1007/978-1-4939-7486-3_13}{\texttt{DOI}}.

\bibitem{Faugeras-version-1}
{\sc O.~Faugeras, L.~Quan, and P.~Sturm}, {\em \titleColor{Self-calibration of a 1d projective camera and its application to the self-calibration of a 2d projective camera}}, in Computer Vision --- ECCV'98, H.~Burkhardt and B.~Neumann, eds., Berlin, Heidelberg, 1998, Springer Berlin Heidelberg, pp.~36--52.

\bibitem{Fulton-intersection-theory}
{\sc W.~Fulton}, \titleColor{{\em Intersection theory}}, vol.~2 of Ergebnisse der Mathematik und ihrer Grenzgebiete. 3. Folge. A Series of Modern Surveys in Mathematics [Results in Mathematics and Related Areas. 3rd Series. A Series of Modern Surveys in Mathematics], Springer-Verlag, Berlin, second~ed., 1998.
\newblock \href{https://doi.org/10.1007/978-1-4612-1700-8}{\texttt{DOI}}.

\bibitem{GKZ}
{\sc I.~Gelfand, M.~Kapranov, and A.~Zelevinsky}, \titleColor{{\em Discriminants, Resultants, and Multidimensional Determinants}}, Modern Birkh{\"a}user Classics, Birkh{\"a}user Boston, 2009.
\newblock \href{https://books.google.com/books?id=ZxeQBAAAQBAJ}{\texttt{URL}}.

\bibitem{Goldman-pyramid-algorithms}
{\sc R.~Goldman}, \titleColor{{\em Pyramid Algorithms}}, The Morgan Kaufmann Series in Computer Graphics, Morgan Kaufmann, San Francisco, 2003.
\newblock \href{https://doi.org/https://doi.org/10.1016/B978-1-55860-354-7.50012-X}{\texttt{DOI}}.

\bibitem{M2}
{\sc D.~R. Grayson and M.~E. Stillman}, {\em \titleColor{Macaulay2, a software system for research in algebraic geometry}}.
\newblock Available at \url{http://www2.macaulay2.com}.

\bibitem{Harris03062018}
{\sc C.~Harris and D.~Lowengrub}, {\em \titleColor{The {C}hern-{M}ather class of the multiview variety}}, Comm. Algebra, 46 (2018), pp.~2488--2499.
\newblock \href{https://doi.org/10.1080/00927872.2017.1392545}{\texttt{DOI}}.

\bibitem{HartleySchaffalitzky}
{\sc R.~Hartley and F.~Schaffalitzky}, {\em \titleColor{Reconstruction from projections using grassmann tensors}}, International Journal of Computer Vision, 83 (2009), pp.~274--293.
\newblock \href{https://doi.org/10.1007/s11263-009-0225-1}{\texttt{DOI}}.

\bibitem{Hartley-Zisserman-book-multiple-view-geometry}
{\sc R.~Hartley and A.~Zisserman}, \titleColor{{\em Multiple view geometry in computer vision}}, Cambridge University Press, Cambridge, second~ed., 2003.
\newblock With a foreword by Olivier Faugeras, \href{https://doi.org/https://doi.org/10.1017/CBO9780511811685}{\texttt{DOI}}.

\bibitem{Ito_Miura_Ueda_2020}
{\sc A.~Ito, M.~Miura, and K.~Ueda}, {\em \titleColor{Projective reconstruction in algebraic vision}}, Canadian Mathematical Bulletin, 63 (2020), p.~592–609.
\newblock \href{https://doi.org/10.4153/S0008439519000687}{\texttt{DOI}}.

\bibitem{MR3705779}
{\sc J.~Kileel}, {\em \titleColor{Minimal problems for the calibrated trifocal variety}}, SIAM J. Appl. Algebra Geom., 1 (2017), pp.~575--598.
\newblock \href{https://doi.org/10.1137/16M1104482}{\texttt{DOI}}.

\bibitem{KL1972-Schubert}
{\sc S.~L. Kleiman and D.~Laksov}, {\em \titleColor{Schubert calculus}}, The American Mathematical Monthly, 79 (1972), pp.~1061--1082.
\newblock \href{http://www.jstor.org/stable/2317421}{\texttt{URL}}.

\bibitem{Lai-Lim-Ye-2025}
{\sc Z.~Lai, L.-H. Lim, and K.~Ye}, {\em \titleColor{Euclidean distance degree in manifold optimization}}, SIAM Journal on Optimization, 35 (2025), pp.~2402--2422.
\newblock \href{https://doi.org/10.1137/25M1735032}{\texttt{DOI}}.

\bibitem{LLY2025-gro}
{\sc Z.~Lai, L.-H. Lim, and K.~Ye}, {\em \titleColor{Grassmannian optimization is {NP}-hard}}, SIAM J. Optim., 35 (2025), pp.~1939--1962.
\newblock \href{https://doi.org/10.1137/24M1672596}{\texttt{DOI}}.

\bibitem{MLB1988-algebra}
{\sc S.~Mac~Lane and G.~Birkhoff}, \titleColor{{\em Algebra}}, Chelsea Publishing Co., New York, third~ed., 1988.
\newblock \href{https://bookstore.ams.org/chel-330/}{\texttt{URL}}.

\bibitem{MR3636592}
{\sc E.~V. Martyushev}, {\em \titleColor{On some properties of calibrated trifocal tensors}}, J. Math. Imaging Vision, 58 (2017), pp.~321--332.
\newblock \href{https://doi.org/10.1007/s10851-017-0712-x}{\texttt{DOI}}.

\bibitem{MRW-multiview}
{\sc L.~G. Maxim, J.~I. Rodriguez, and B.~Wang}, {\em \titleColor{Euclidean distance degree of the multiview variety}}, SIAM J. Appl. Algebra Geom., 4 (2020), pp.~28--48.
\newblock \href{https://doi.org/10.1137/18M1233406}{\texttt{DOI}}.

\bibitem{MRW2021-edprojective}
{\sc L.~G. Maxim, J.~I. Rodriguez, and B.~Wang}, {\em \titleColor{Euclidean distance degree of projective varieties}}, Int. Math. Res. Not. IMRN,  (2021), pp.~15788--15802.
\newblock \href{https://doi.org/10.1093/imrn/rnz266}{\texttt{DOI}}.

\bibitem{MRW-survey}
{\sc L.~t.~G. Maxim, J.~I. Rodriguez, and B.~Wang}, {\em \titleColor{Applications of singularity theory in applied algebraic geometry and algebraic statistics}}, in Handbook of geometry and topology of singularities {VII}, Springer, Cham, 2025, pp.~767--818.
\newblock \href{https://doi.org/10.1007/978-3-031-68711-2\_14}{\texttt{DOI}}.

\bibitem{Miller-Sturmfels-combinatorial}
{\sc E.~Miller and B.~Sturmfels}, \titleColor{{\em Combinatorial commutative algebra}}, vol.~227 of Graduate Texts in Mathematics, Springer-Verlag, New York, 2005.
\newblock \href{https://doi.org/https://doi.org/10.1007/b138602}{\texttt{DOI}}.

\bibitem{LongKanade}
{\sc L.~Quan and T.~Kanade}, {\em \titleColor{Affine structure from line correspondences with uncalibrated affine cameras}}, IEEE Transactions on Pattern Analysis and Machine Intelligence, 19 (1997), pp.~834--845.
\newblock \href{https://doi.org/10.1109/34.608285}{\texttt{DOI}}.

\bibitem{rydell2024projectionshigherdimensionalsubspaces}
{\sc F.~Rydell}, {\em \titleColor{Projections of higher dimensional subspaces and generalized multiview varieties}}, 2024.
\newblock \href{https://arxiv.org/abs/2309.10262}{\texttt{arXiv:} \texttt{2309.10262}}.

\bibitem{Rydell_2023_ICCV}
{\sc F.~Rydell, E.~Shehu, and A.~Torres}, {\em \titleColor{Theoretical and numerical analysis of 3d reconstruction using point and line incidences}}, in Proceedings of the IEEE/CVF International Conference on Computer Vision (ICCV), October 2023, pp.~3748--3757.
\newblock \href{https://openaccess.thecvf.com/content/ICCV2023/html/Rydell_Theoretical_and_Numerical_Analysis_of_3D_Reconstruction_Using_Point_and_ICCV_2023_paper.html}{\texttt{URL}}.

\bibitem{SSN2005-three-view}
{\sc H.~Stewenius, F.~Schaffalitzky, and D.~Nister}, {\em \titleColor{How hard is 3-view triangulation really?}}, in Tenth IEEE International Conference on Computer Vision (ICCV'05) Volume 1, vol.~1, 2005, pp.~686--693 Vol. 1.
\newblock \href{https://doi.org/10.1109/ICCV.2005.115}{\texttt{DOI}}.

\bibitem{Wolf}
{\sc L.~Wolf and A.~Shashua}, {\em \titleColor{On projection matrices and their applications in computer vision}}, in IEEE International Conference on Computer Vision, 2001.
\newblock \href{https://doi.ieeecomputersociety.org/10.1109/ICCV.2001.10057}{\texttt{URL}}.

\end{thebibliography}

\bigskip \medskip \bigskip

\noindent
\footnotesize {\bf Authors' addresses:}
\smallskip

\noindent Bella Finkel, University of Wisconsin--Madison, USA \hfill {\tt  blfinkel@wisc.edu} \url{https://people.math.wisc.edu/~blfinkel/}

\smallskip
\noindent Jose Israel Rodriguez, University of Wisconsin--Madison, USA \hfill {\tt  jose@math.wisc.edu}\newline
\url{https://sites.google.com/wisc.edu/jose/}
\end{document}